\documentclass[preprint,3p,12pt,pdf]{elsarticle}

\usepackage{lineno,hyperref}
\modulolinenumbers[5]

\usepackage{hyperref}

\usepackage{epstopdf}

\usepackage{tikz}
\usetikzlibrary{arrows}

\usepackage{pst-node}
\usepackage{tikz-cd}

\usepackage[shellescape]{gmp}

\usepackage{mathrsfs}
\usepackage{amssymb}
\usepackage{amsfonts}
\usepackage{latexsym}
\usepackage{mathtools}
\usepackage{xcolor}
\usepackage{wrapfig}
\usepackage{floatflt}

\usepackage{mathtools}
\usepackage{extarrows}

\usepackage{graphicx}

\usepackage{subcaption}

\usepackage{endnotes}

\def\lb{\label}

\newcommand{\er}[1]{\textrm{(\ref{#1})}}

\newtheorem{theorem}{\bf Theorem}[section]

\def\a{\alpha}

\def\d{\delta}         
\def\D{\Delta}         
           
\def\z{\zeta}   \def\cH{{\mathcal H}}       \def\mH{{\mathscr H}}
           
\def\p{\psi}

\def\m{\mu}            
            
\def\r{\rho}

\def\ve{\varepsilon}       \def\vp{\varphi}    

\def\Z{{\mathbb Z}}    \def\R{{\mathbb R}}       
    \def\N{{\mathbb N}}   

\def\lt{\biggl}                  \def\rt{\biggr}
\def\ol{\overline}

\let\ge\geqslant                 \let\le\leqslant

\def\iy{\infty}

                 \def\ev{\equiv}
        
\def\el2{\ell^{\,2}}             \def\1{1\!\!1}

\def\const{\mathop{\mathrm{const}}\nolimits}

\def\Re{\mathop{\mathrm{Re}}\nolimits}

\def\BBox{\hspace{1mm}\vrule height6pt width5.5pt depth0pt \hspace{6pt}}

\newtheorem{corollary}[theorem]{\bf Corollary}

\let\ge\geqslant
\let\le\leqslant

\newcommand{\ca}{\begin{cases}}
	\newcommand{\ac}{\end{cases}}
\newcommand{\ma}{\begin{pmatrix}}
	\newcommand{\am}{\end{pmatrix}}
\def\eq{\begin{equation}}
	\def\qe{\end{equation}}

\renewcommand{\[}{\begin{equation}}
	\renewcommand{\]}{\end{equation}}

\def\BBox{\hspace{1mm}\vrule height6pt width5.5pt depth0pt \hspace{6pt}}

\begin{document}
	
	\begin{frontmatter}

		\title{Generalized Schr\"oder-type functional equations for Galton--Watson processes in random environments}

		\date{\today}

		\author
		{Anton A. Kutsenko}
	
		\address{Mathematical Institute for Machine Learning and Data
Science, KU Eichst\"att--Ingolstadt, Germany; email: akucenko@gmail.com}
	
	\begin{abstract}
		  The classical Galton--Watson process works with a fixed probability of fission at each time step. One of the generalizations is that the probabilities depend on time. We consider one of the most complex and interesting cases when we do not know the exact probabilities of fission at each time step - these probabilities are random variables themselves. The limit distributions of the number of descendants are described in terms of generalized integral and differential functional equations of the Schr\"oder type. There are no more analogs of periodic Karlin-McGregor functions, which were very helpful in the analysis of the asymptotic behavior of limit distributions for the classical case. We propose some approximate asymptotic methods. Even simple cases of random families with one or two members lead to nice asymptotics involving, interesting problems related to special functions and special constants. One of them, Example 2 is already announced on \href{https://math.stackexchange.com/questions/4748129/asymptotics-of-sequence-of-rational-numbers}{this}\endnote{\it https://math.stackexchange.com/questions/4748129/asymptotics-of-sequence-of-rational-numbers\lb{ref1}} and \href{https://mathoverflow.net/questions/458885/simple-integral-equation}{this}\endnote{\it https://mathoverflow.net/questions/458885/simple-integral-equation\lb{ref2}} sites. Finally, the phenomenon of why the oscillations in the main asymptotic term usual for the classical case become rare in the case of random environments is discussed.
	\end{abstract}


	
\end{frontmatter}


{\section{Introduction}\lb{sec0}}

A simple Galton--Watson process, see, e.g. \cite{H}, is defined by
\[\lb{001}
 X_{t+1}=\sum_{j=1}^{X_t}\xi_{j,t},\ \ \ X_0=1,\ \ \ t\in\N\cup\{0\},
\]
where all $\xi_{j,t}$ are independent and identically-distributed natural number-valued random variables with the probability-generating function
\[\lb{002}
 P(z):=\mathbb{E}z^{\xi}=p_0+p_1z+p_2z^2+p_3z^3+....
\]
All $p_j$ are non-negative and their sum equals $P(1)=1$. We focus on the so-called supercritical case $p_0=0$ and $0<p_1<1$. In this case, one may define the limit distribution
\[\lb{003}
 \vp_n:=\lim_{t\to+\iy}\frac{\mathbb{P}(X_t=n)}{\mathbb{P}(X_t=1)}
\] 
describing the so-called relative limit densities $\vp_n$, see \cite{H}. These values are coefficients of the analytic function
\[\lb{004}
 \Phi(z)=\vp_1z+\vp_2z^2+\vp_3z^3+....
\]
The function $\Phi$ satisfies the Scr\"oder-type functional equation
\[\lb{005}
 \Phi(P(z))=p_1\Phi(z),\ \ \ \Phi(0)=0,\ \ \ \Phi'(0)=1,
\]
it is defined at least for $|z|<1$.
A Galton--Watson process in a random environment is characterized by the fact that $\xi_{j,t}$ are still independent but they are identically distributed for fixed $t$ only. The simplest case is when $P(z)$, see \er{002}, is fixed for each time $t$, but it depends on $t$ deterministically, see, e.g., \cite{J}. We consider a more complex case when the dependence on $t$ is random. This means that $P(z)$ is again not fixed and, at each time step $t$, it is chosen randomly among
\[\lb{006}
 P_r(z)=p_{1r}z+p_{2r}z^2+p_{3r}z^3+...,
\] 
where $r\in R$ is a random parameter having a distribution measure $\m$ defined on some probabilistic space $R$. Note that some other generalizations, e.g., when $(p_n(t))_{n=0}^{+\iy}$ are independent identically distributed random vectors, are considered in, e.g., \cite{VD}, see also references therein. For the case \er{006}, we can still define the limit distribution \er{003} and the corresponding analytic function \er{004}. The difference between simple and random cases is that $\Phi$ satisfies the generalized Schr\"oder-type functional equation.

\begin{theorem}\lb{T1}
Let $\Phi(z)$ be constructed as an analog of \er{003} and \er{004} for the random environments with generating functions \er{006}. If $\int_Rp_{1r}\m(dr)>0$ and $\frac{\int_Rp_{1r}^2\m(dr)}{\int_Rp_{1r}\m(dr)}<1$ then $\Phi(z)$ satisfies the following integro-functional equation
\[\lb{007}
 \int_R\Phi(P_r(z))\m(dr)=\Phi(z)\int_Rp_{1r}\m(dr),\ \ \ \Phi(0)=0,\ \ \ \Phi'(0)=1.
\]
Moreover, $\Phi(z)$ is analytic at least for small $z$. 
\end{theorem}
{\bf Proof.} Consider the complete metric space of analytic functions defined for $|z|\le h$ and equipped with the standard max-norm,
\[\lb{008b}
 \mH_{h,C}=\{f(z):\ |f(z)-z|\le C|z|^2,\ for\ |z|\le h\}.
\]
Consider the following operator
\[\lb{008c}
 \cH(f)(z)=\frac{\int_Rf(P_r(z))\m(dr)}{\int_Rp_{1r}\m(dr)}.
\]
It is seen that $|P_r(z)|\le P_r(|z|)\le|z|$ for $|z|\le1$, since $P_r(0)=0$, $P_r(1)=1$, and all their Taylor coefficients are non-negative. Thus, 
for any $0<h<1$, operator $\cH$ is defined on $\mH_{h,C}$. Moreover,
\[\lb{008d}
 \cH(f)(z)-z=\frac{\int_R(f(P_r(z))-p_{1r}z)\m(dr)}{\int_Rp_{1r}\m(dr)},
\]
which leads to
\[\lb{008e}
|\cH(f)(z)-z|\le \frac{\int_R(|P_r(z)-p_{1r}z|+C|P_r(z)|^2)\m(dr)}{\int_Rp_{1r}\m(dr)}\le A|z|^2+C\frac{\int_Rp_{1r}^2\m(dr)}{\int_Rp_{1r}\m(dr)}|z|^2(1+B|z|).
\]
for some $A,B>0$ and $|z|\le h<1$, since $P_{r}(z)=p_{1r}z+Q_r(z)$ for some $Q_r(z)$ satisfying $|Q_r(z)|\le|z|^2$ for $|z|\le1$. Recall that all the coefficients of $P_r$ are non-negative and their sum is equal to $1$. Condition $\frac{\int_Rp_{1r}^2\m(dr)}{\int_Rp_{1r}\m(dr)}<1$ allows us to take sufficiently large $C>0$ and small $h>0$ such that $\cH(f)\in\mH_{h,C}$ for $f\in\mH_{h,C}$. Using the same arguments as above, it is easy to see that
\[\lb{008g}
|\cH(f)(z)-\cH(g)(z)|\le B|z|^2,
\]
where, again, $B>0$ is some constant independent on the choice of functions $f,g\in\mH_{h,C}$. Thus, $\cH$ is a contraction mapping if $h>0$ is small enough. Hence, there is a unique fixed point $\Phi$ - the solution of \er{007}. Now, we have
\begin{multline}\lb{008a}
 \Phi_t(z):=\frac{\mathbb{E}z^{X_t}}{\mathbb{P}(X_t=1)}=\frac{\int_{R^t}P_{r_t}\circ...\circ P_{r_1}(z)\m(dr_1)...\m(dr_t)}{\int_{R^t}p_{1r_t}...p_{1r_1}\m(dr_1)...\m(dr_t)}=\\
 \frac{\int_{R^t}P_{r_t}\circ...\circ P_{r_1}(z)\m(dr_1)...\m(dr_t)}{(\int_{R}p_{1r}\m(dr))^t}=\frac{\int_R\Phi_{t-1}(P_{r_1}(z))\m(dr_1)}{\int_Rp_{1r}\m(dr)}.
\end{multline}
Taking into account the fact that $\Phi_0(z)=z\in\mH_{h,C}$, we obtain that $\Phi_t(z)$ converges to the fixed point of $\cH$ as standard sequential iterations in the Banach Fixed Point theorem.  \BBox

Usually, $\Phi$ can be analytically extended to, at least, the open unit disk.

\begin{corollary}\lb{C1} If $\int_Rp_{1r}\m(dr)>0$ and $p_{1r}$ is uniformly separated from $1$ then $\Phi(z)$ can be analytically extended to the open unit disk. Moreover, there is $\a\in\R$ such that $|\Phi(z)|<(1-|z|)^{\a}$.
\end{corollary}
{\bf Proof.} There is $0<\d<1$ such that $p_{1r}<1-\d$. All the conditions of Theorem \ref{T1} are satisfied. If $|z|\le1$ then
\[\lb{009a}
 |P_r(z)|\le p_{1r}|z|+(1-p_{1r})|z|^2\le(1-\d)|z|+\d|z|^2.
\]
The assumption
\[
 |z|<Rh,\ \ \ R:=\frac{2}{1-\d+\sqrt{(1-\d)^2+4h\d}}
\]
leads to
\[\lb{009c}
 (1-\d)|z|+\d|z|^2<h,
\]
and, hence, $|P_r(z)|<h$, see \er{009a}. Thus, if $\Phi(z)$ is analytic for $|z|<h<1$ then it is also analytic for $|z|<Rh$ by \er{007}. Hence,  we can increase the domain of analyticity, since $R>1$ for $h<1$. Using, e.g., \er{008a}, we see that all the Taylor coefficients $\vp_n$ of $\Phi(z)$ are non-negative. Thus, it is sufficient to prove $\Phi(x)<(1-x)^{\a}$ for $x\in[0,1]$. For small $x\le h$, the inequality is obvious for any $\a<0$. Using the same arguments as above, we can increase the interval in $R>1$ times following \er{007}: 
\[\lb{009d}
\Phi(x)=\frac{\int_R\Phi(P_r(x))\m(dr)}{\int_Rp_{1r}\m(dr)}<(1-x)^{\a}\frac{\int_R\lt(\frac{1-P_r(x)}{1-x}\rt)^{\a}\m(dr)}{\int_Rp_{1r}\m(dr)}\le(1-x)^{\a},
\]
where $\a<0$ is chosen such that the following inequality is fulfilled 
\begin{multline}\lb{009e}
 \lt(\frac{1-P_r(x)}{1-x}\rt)^{\a}=(p_{1r}+p_{2r}(1+x)+p_{3r}(1+x+x^2)+...)^{\a}
 =(1+p_{2r}x+p_{3r}(x+x^2)+...)^{\a}=\\
 (1+(1-p_{1r})x+p_{3r}x^2+p_{4r}(x^2+x^3)...)^{\a}\le(1+\d x)^{\a}<\int_Rp_{1r}\m(dr)
\end{multline}
for, say, $x\in[h/2,h]$. It is possible, since $\frac{1-P_r(x)}{1-x}$ is strictly greater than $1$ for $x\in[h/2,h]$, and we can take, e.g., $\a=\min\{-1,\log_{1+\d h/2}\int_Rp_{1r}\m(dr)\}$, see \er{009e}. Next interval increasing steps can be performed with the same $\a<0$, since $\frac{1-P_r(x)}{1-x}$ is a monotonically increasing function. \BBox

The condition on separation of $p_{1r}$ from $1$ can be weakened in many cases with the help of taking limits $p_{1r}\to1$ for required $r$, plus a little modification of other $p_{jr}$. Also, further improvements may allow us to extend the domain of definition of $\Phi$ to the intersection of filled Julia sets related to $P_r$ - this is seen from \er{008a}. 

Differentiating \er{007} at $z=0$ and using Fa\`a di Bruno's formula for the composition of two functions, we derive explicit recurrence expressions for the Tailor coefficients $\vp_n$ of $\Phi(z)$, see \er{003} and \er{004},
\[\lb{008}
 \vp_n=\lt(\int_R(p_{1r}-p_{1r}^{n})\m(dr)\rt)^{-1}\sum_{k=1}^{n-1}\vp_{k}\int_RB_{n,k}(p_{1r},2p_{2r},...,(n-k+1)!p_{n-k+1,r})\m(dr),
\]
where $B_{n,k}$ are Bell polynomials
\[\lb{Bell1}
B_{n,k}(x_1,x_2,...,x_{n-k+1}):=\sum\frac{n!}{j_1!j_2!...j_{n-k+1}!}\lt(\frac{x_1}{1!}\rt)^{j_1}\lt(\frac{x_2}{2!}\rt)^{j_2}...\lt(\frac{x_{n-k+1}}{(n-k+1)!}\rt)^{j_{n-k+1}},
\]
and the sum for Bell polynomials is taken over all sequences $j_r$ of non-negative integers such that
\[\lb{Bell2}
j_1+j_2+...+j_{n-k+1}=k,\ \ \ j_1+2j_2+3j_3+...+(n-k+1)j_{n-k+1}=n.
\] 
The natural assumption is that all the integrals exist in \er{008}. 

\begin{figure}
	\centering
	\includegraphics[width=0.6\linewidth]{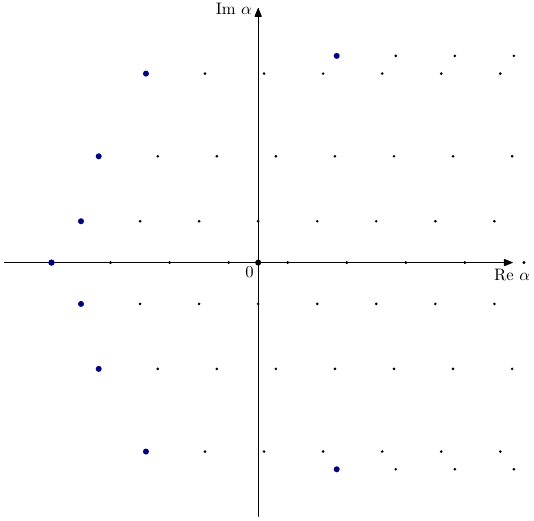}
	\caption{Example of distribution of $\a$ in \er{019}. Large blue points are zeros of \er{015}.}
	\label{figr}
\end{figure}

In some cases, especially when $P_r$ depends on $r$ linearly, it is convenient to rewrite \er{007} in a differential form introducing $F(z)=\int_0^z\Phi(\z)d\z$. We consider this case in the Example section. Integral and differential extensions of the Schr\"oder-type functional equations are much more complex than the original form, where $P_r$ does not depend on $r$. In particular, it is hard to determine the asymptotic behavior of $\vp_n$ for large $n\to+\iy$. One of the reasons is that there are no more direct analogs of Karlin--McGregor functions defined in \cite{KM1} and \cite{KM2} for the classical case. Let us recall some facts about this function. For the case $P_r\ev P$, it is possible to define
\[\lb{009}
 \Pi(z):=\lim_{t\to+\iy}\underbrace{P\circ...\circ P}_{t}(1-E^{-t}z),
\]
where $E=P'(1)$ is the expectation $E=p_1+2p_2+3p_3+...$. The function $\Pi(z)$ is entire when $P$ is entire, and it satisfies the Poincar\'e-type functional equation
\[\lb{010}
 P(\Pi(z))=\Pi(Ez),\ \ \ \Pi(0)=1,\ \ \ \Pi'(0)=-1.
\]
Note that this function is related to the probability density function for the so-called {\it martingal limit}, which is another limit distribution for branching processes apart from \er{003}:
\[\lb{011}
 p(x):=\lim_{t\to+\iy}E^t\mathbb{P}(X_t=[xE^t])=\frac1{2\pi\mathbf{i}}\int_{\mathbf{i}\R}\Pi(z)e^{zx}dz,
\]
where the square brackets denote the integer part of a number, see, e.g. \cite{D1}. Combining $\Phi$ and $\Pi$, and using their functional equations, we can define the $1$-periodic function
\[\lb{012}
 K(z):=\Phi(\Pi(E^z))p_1^{-z}.
\]
This function is very efficient in analyzing asymptotic series for both limit distributions: \er{003} and \er{011}.  The direct analog of \er{009} is
\[\lb{013}
 \Pi(z):=\lim_{t\to+\iy}\int_{R^t}\underbrace{P_{r_1}\circ...\circ P_{r_t}}_{t}(1-E^{-t}z)\m(dr_t)...\m(dr_1),\ \ \ \Pi(0)=1,\ \ \ \Pi'(0)=-1,
\]
where the average expectation $E$ is given by
\[\lb{014}
 E:=\int_R(P_r)'(1)\m(dr).
\]
However, the definition of \er{013} is already questionable - for the case  $(P_r)'(1)\ne \const$, try to compute the second derivative of $\Pi$.   The analog of \er{011} is also questionable in this case. Another variant is something like
\[\lb{013a}
 \Pi(z):=\lim_{t\to+\iy}\int_{R^t}\underbrace{P_{r_1}\circ...\circ P_{r_t}}_{t}\lt(1-\frac{z}{\prod_{j=1}^t(P_{r_j})'(1)}\rt)\m(dr_t)...\m(dr_1),\ \ \ \Pi(0)=1,\ \ \ \Pi'(0)=-1,
\]
but we do not plan to discuss the details. In both cases \er{013} and \er{013a}, there are no appropriate variants of \er{010}, making the definition of $1$-periodic analogs of Karlin--McGregor functions \er{012} almost impossible. This fact significantly complicates the asymptotic analysis of the coefficients $\vp_n$ when $n\to+\iy$. However, a few methods allow us to obtain good approximations of $\vp_n$. Often, the main contribution to the asymptotic behavior of $\vp_n$ comes from the singularity of $\Phi(z)$ at $z=1$. Substituting {\it ansatz} $\Phi(z)\approx A(1-z)^{\a}$ into the main equation \er{007}, we obtain the equation for $\a$:
\[\lb{015}
 \int_{R}((P_r)'(1))^{\a}\m(dr)=\int_Rp_{1r}\m(dr).
\] 
Usually, \er{015} has infinitely many solutions $\a$, which are zeros of the corresponding (entire) function. The largest contribution comes from $\a$ with the smallest real parts, because
\[\lb{016}
 \Phi(z)\approx A(1-z)^{\a}=A\sum_{n=0}^{+\iy}(-1)^n\binom{\a}{n} z^n,
\]
and using the asymptotic for generalized binomial coefficients
\[\lb{017}
\binom{\a}{n}\simeq\frac{(-1)^n}{\Gamma(-\a)n^{\a+1}}+\frac{(-1)^n\a(\a+1)}{2\Gamma(-\a)n^{\a+2}}+...,
\]
we obtain the final approximation
\[\lb{018}
 \vp_n\approx \sum C_{\a}n^{-\a-1},
\]
where the sum is taken over the solutions of \er{015} with small real parts. Because $P_r'(1)\ge1$ and $p_{1r}$ are non-negative, the real parts of solutions $\a$ are bounded from below. Note that if $\a=\a_1+\mathbf{i}\a_2$ is complex then the corresponding term in \er{018} gives the long-phase logarithmic oscillation $C_{\a}n^{-\a_1-1}\exp(-\mathbf{i}\a_2\ln n)$ (plus the corresponding terms related to the complex conjugate $\ol{\a}$) similar to that in the asymptotic of the classical Galton--Watson process with one generating function, see, e.g., \cite{K} and \cite{DG}. However, in the mixed case of many generating functions \er{007}, the derivation of complete asymptotic series for power coefficients of $\Phi(z)$ is, generally speaking, an open problem. Summarizing the above, we can formulate the general rough idea of finding approximations: Substitute
\[\lb{019}
 \vp_n\simeq\sum_{\a\in\{Solutions\ of\ \er{015}\},\ j\in\N}\frac{C_{\a,j}}{n^{\a+j}}
\] 
into \er{008} directly or \er{016} into \er{007} and use \er{017} with some {\it a priori} and/or empirical additional information to obtain relations for coefficients $C_{\a,j}$. However, in many cases, it does not work properly and more sophisticated tricks are required. An example of the distribution of $\a$ in \er{019} is illustrated in Fig. \ref{figr}. 

\begin{figure}
	\centering
	\includegraphics[width=0.5\linewidth]{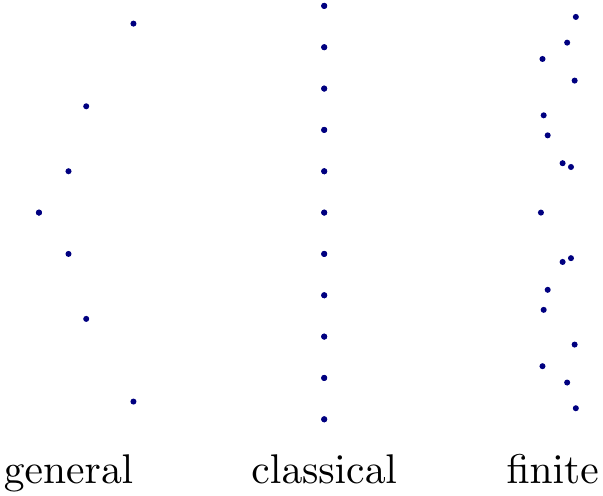}
	\caption{Typical distributions of zeros of \er{015} in three common cases: classical - one polynomial; finite - finite number of polynomials; general - all other cases, e.g. $\m$ is continuous.}
	\label{figra}
\end{figure}

We will treat the explained approximate method in the examples below. Moreover, in some cases, we take more approximate terms for $\Phi(z)$ at $z=1$ and other critical points of the unit circle.
We consider four examples. For the simplest linear fractional generating functions, which correspond to the power-law distribution, \er{007} can be solved explicitly. This is the complete analog of the similar case for the classical Galton--Watson process. The next two examples work with families of one or two particles. At each time step the generating function \er{006} is simply $P_r(z)=rz+(1-r)z^2$, where $r$ is distributed uniformly in $(1/2,1)$ or in $(0,1)$. This means that at each time step each particle generates one new particle or stays alone. The generation probabilities are any from $(1/2,1)$ or $(0,1)$ uniformly. Such types of problems are natural in situations when we do not know the exact probabilities of generation (or fission)  at each time step but expect them with some probabilities. Two cases with distribution intervals $(1/2,1)$ and $(0,1)$ differ significantly. Relative limit densities \er{003} for both cases satisfy power laws with different exponents. In contrast to the classical case, see, e.g., \cite{K}, it seems that there are no long-phase oscillations for the first approximation term. These long-phase oscillations may appear in the next asymptotic terms. However, in contrast to the half-interval $(1/2,1)$, for the full interval $(0,1)$ there are short-phase oscillations with bounded periods. In the last example, we consider the uncertainty produced by only two polynomials. This case is closest to the classical point of one polynomial and shows the long-phase logarithmic oscillations usual for the classical cases. As mentioned in, e.g., \cite{DIL, DMZ, CG}, such short- and long-phase oscillations can be important in applications in physics and biology. Let us discuss oscillations a little bit more. Due to \er{018} and the arguments below, they appear in some asymptotic terms for $\vp_n$. However, for the classical Galton--Watson process, the oscillations appear for the main asymptotic term already, because $\a=(\ln P'(0)+2\pi\mathbf{i}m)/\ln P'(1)$, $m\in\Z$ are the only roots of \er{015}. Each complex $\a$ has the same real part as the primary real $\a$ and contributes to the main asymptotic term with some long-phase oscillations. In the generalized case, when we have more than one generating function $P_r(z)$, there is still one primary real solution $\a$ of \er{015}. The existence of complex solutions $\a$ of \er{015} with the same real part as the primary one requires that there is some real $B$ such that $P_r(1)'=B^{n_r}$ for some integer $n_r$ and almost all $r$. This condition describes some one-dimensional ``algebraic" curve in $\R^R$ with the parameter-argument $B$. The union of all such curves over all integer vectors $(n_r)$ is dense but has a ``zero measure" in most acceptable models. Let us briefly discuss short-phase oscillations, e.g., considered in Example 2 below. They usually come from other critical points of the unit circle, not equal to $z=1$. However, since the real part of $P_r(z_0)$ with complex $z_0$ is less than $P_r(1)$, see \er{002}, the real part of the corresponding $\a$ should be larger than the minimal real solution $\a$ of \er{015}. Roughly speaking, this means that, for most of the cases, short-phase oscillations will appear not in the main asymptotic term. Thus, oscillations in the main term of the asymptotics of the number of descendants, which is mostly common in the classical branching processes with fixed offspring distributions, become rare for generalized cases taken in a random environment. Two examples related to this discussion are considered at the beginning (Example 1) and end (Example 3) of the paper. To summarize the above observations, let us look at Fig. \ref{figra} with typical distributions of zeros of \er{015} in various cases. In the classical case, the zeros have the same real part - thus, the oscillations appear already in the main asymptotic term. For the case of mixing of finite number of polynomials, typically there are infinitely many $\a$ whose real parts are arbitrarily close to the minimal one. The corresponding oscillations are generally visible in the main asymptotic term, but not so strong as in the classical case. Of course, in some special cases, we can still have infinitely many $\a$ whose real part coincides exactly with the minimal one. In this case, the difference with the classical case of one polynomial is not large. For the general case of, e.g. continuous measure $\m$, we have one (minimal) real zero. The main asymptotic term satisfies a pure power law without any oscillations. 

\section{Examples}

\subsection{Example 0.} In rare cases, \er{007} can be solved explicitly. As in the classical case of the Galton--Watson process, if we assume $P_r(z)=(1-r)z/(1-rz)$ but with a varying distribution of $r$ with a measure $\m(dr)$ then $\Phi(z)=z/(1-z)$ is the solution of \er{007}.

\subsection{Example 1.} 
\begin{figure}
	\centering
	\begin{subfigure}{0.99\linewidth}
		\includegraphics[width=\linewidth]{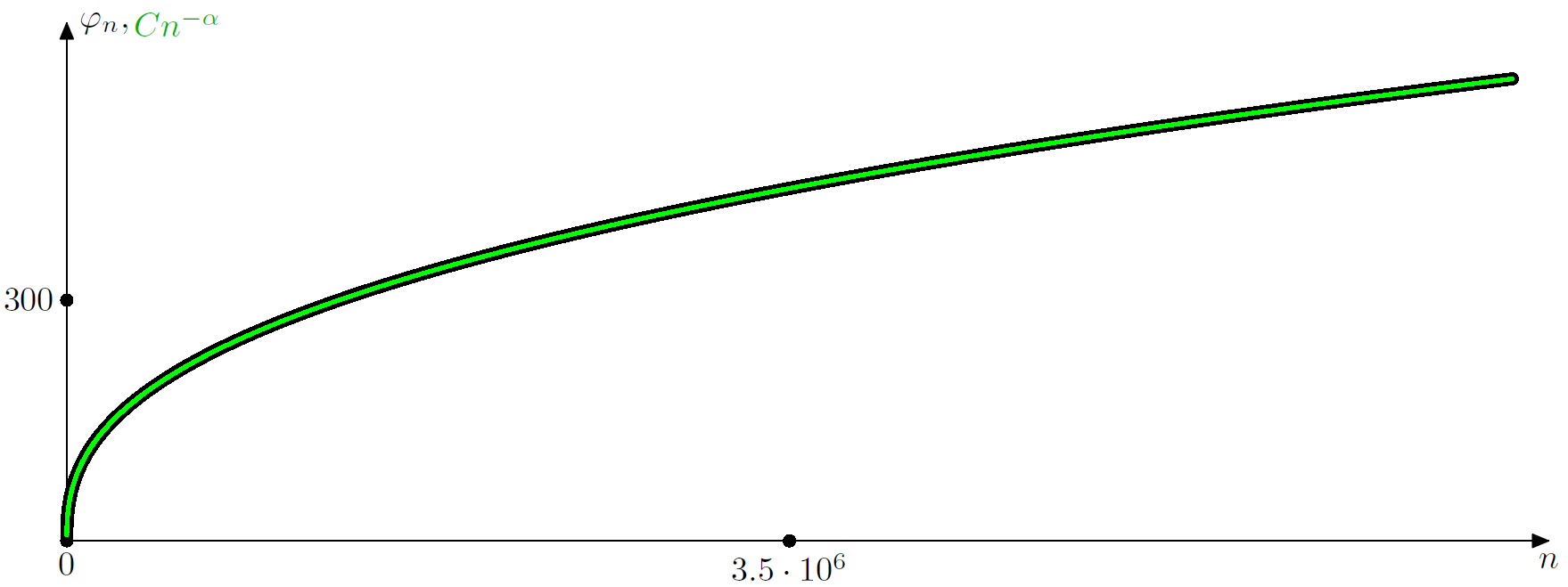}
		\caption{Exact value and first approximation term}
		\label{fig1a}
	\end{subfigure}\vfill
	\begin{subfigure}{0.99\linewidth}
		\includegraphics[width=\linewidth]{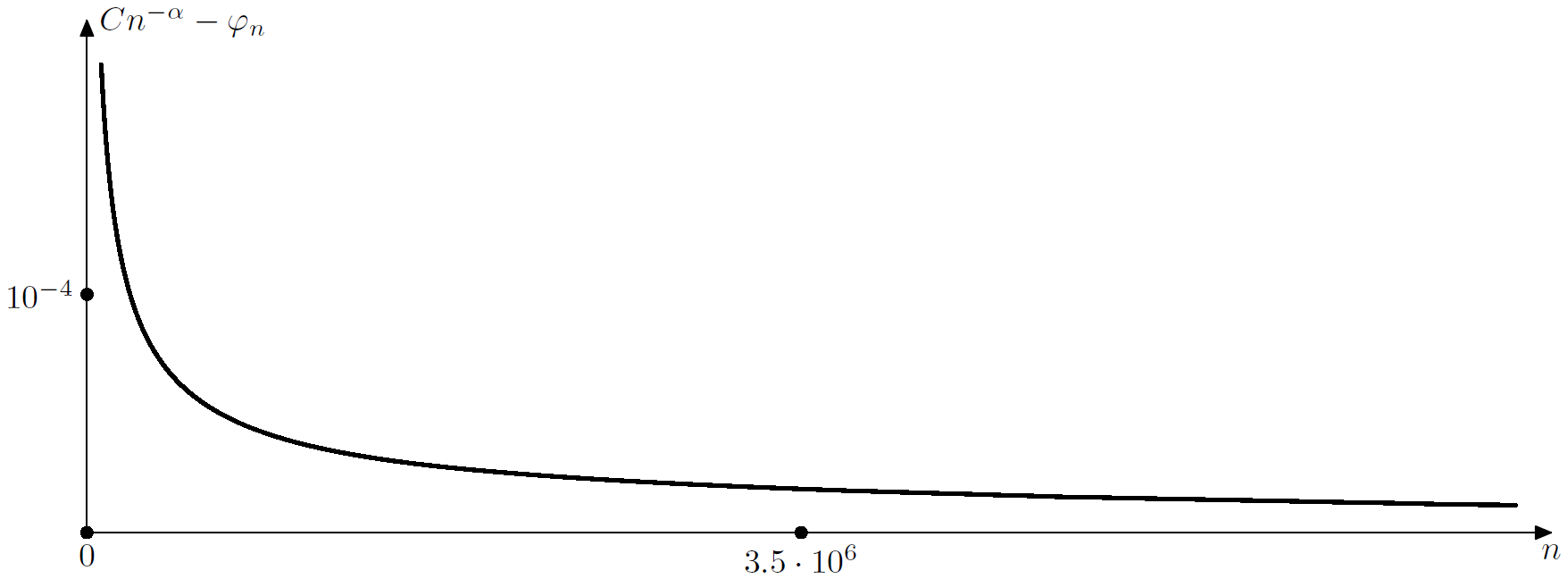}
		\caption{Difference between first approximation term and exact value}
		\label{fig1b}
	\end{subfigure}\vfill
	\begin{subfigure}{0.99\linewidth}
		\includegraphics[width=\linewidth]{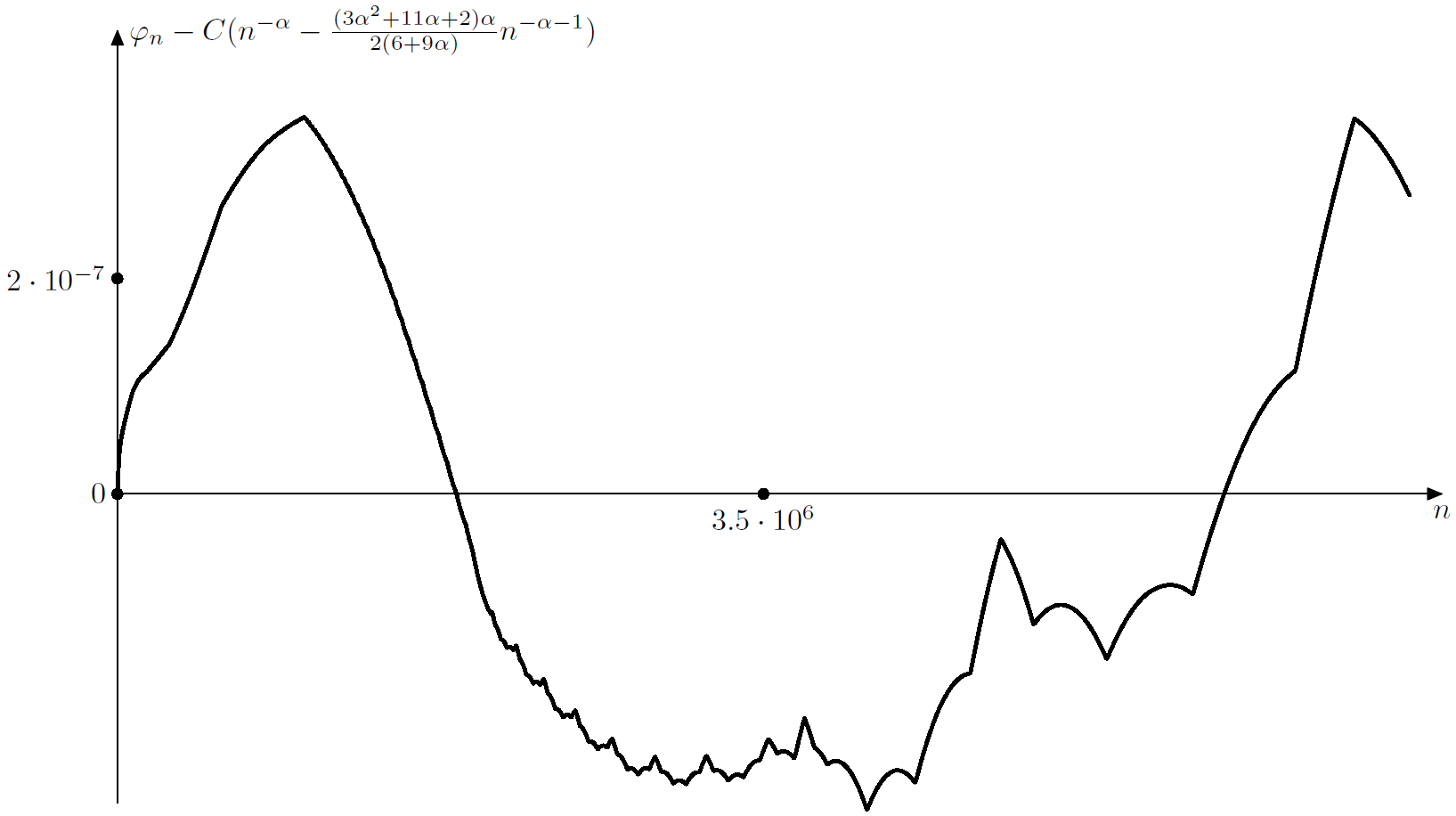}
		\caption{Difference between exact value and two first approximation terms}
		\label{fig1c}
	\end{subfigure}
	\caption{Comparison between exact values and their approximations, see \er{111}.}
	\label{fig1}
\end{figure}

\begin{figure}
	\centering
	\includegraphics[width=0.99\linewidth]{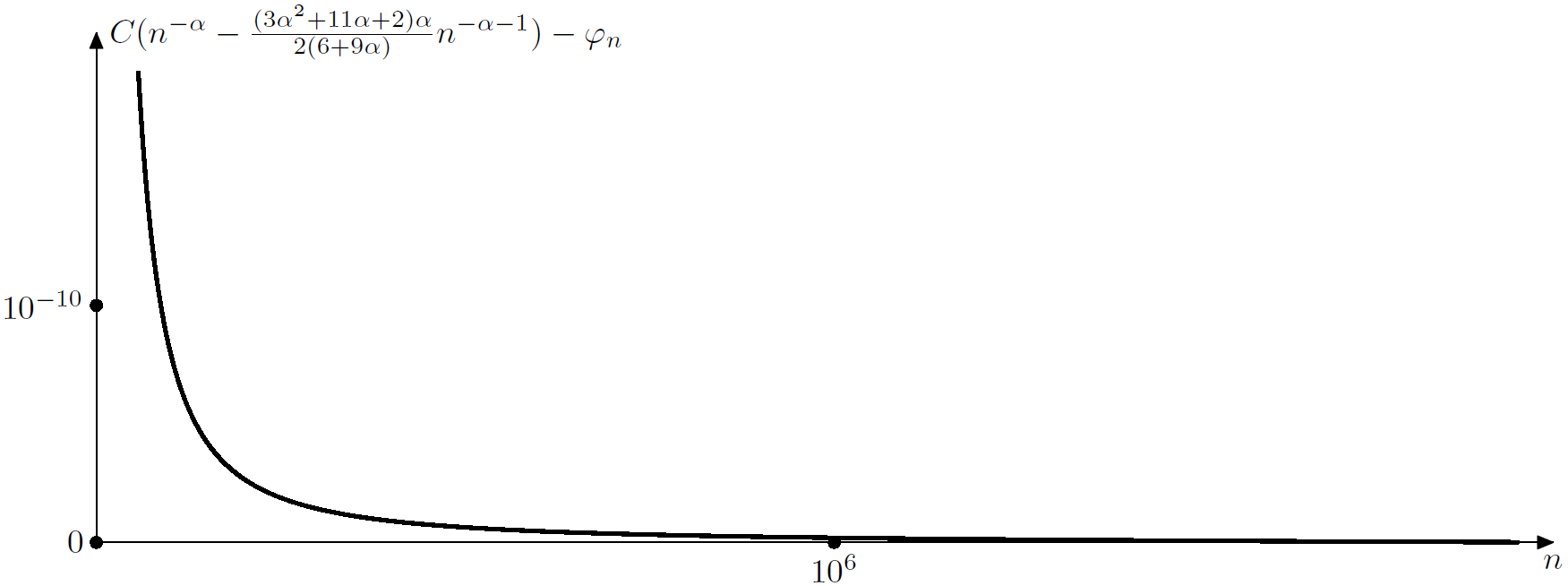}
	\caption{Correction of Fig. \ref{fig1c} with double-double precision and more precise $C$, see \er{110dd}.}
	\label{fig1cdd}
\end{figure}

The case $P_r(z)=rz+(1-r)z^2$, with uniform $r\in(1/2,1)$. Equation \er{007} is
\[\lb{100}
 \int_{\frac12}^{1}\Phi(rz+(1-r)z^2)dr=\frac38\Phi(z).
\]
Taking $F(z)=\int_0^z\Phi(\z)d\z$, \er{100} can be rewritten as
\[\lb{101}
 \frac{F(z)-F(\frac{z+z^2}2)}{z-z^2}=\frac38F'(z),\ \ F(0)=F'(0)=0,\ \ F''(0)=1.
\]
Tailor coefficients, see \er{004}, can be determined from \er{101} after some transformations:
\[\lb{102}
 \sum_{n=1}^{+\iy}\frac{\vp_nz^n}{2(n+1)}\sum_{j=0}^n\lt(\frac{z+z^2}2\rt)^j=\frac38\sum_{n=1}^{+\iy}\vp_nz^n,\ \ \ \vp_1=1,
\]
which leads to
\[\lb{103}
 \frac38\vp_n=\frac{\vp_n}{2(n+1)}\lt(\frac{1}{2^0}\binom{0}{0}+...+\frac1{2^n}\binom{n}{0}\rt)+\frac{\vp_{n-1}}{2n}\lt(\frac1{2^1}\binom{1}{1}+...+\frac1{2^{n-1}}\binom{n-1}{1}\rt)+....
\]
Identity \er{103} gives the recurrence relation
\[\lb{104}
 \lt(\frac38-\frac{1-2^{-n-1}}{n+1}\rt)\vp_n=\frac{\vp_{n-1}}{2n}c_{n,1}+\frac{\vp_{n-2}}{2(n-1)}c_{n,2}+...+\frac{\vp_{n-k}}{2(n-k+1)}c_{n,k},\ \ \ \vp_1=1,
\]
where $n\ge2$ and $k$ is the maximal number such that $2k\le n$. The coefficients $c_{n,j}$ are given by
\[\lb{105}
 c_{n,j}=\frac1{2^{j}}\binom{j}{j}+...+\frac1{2^{n-j}}\binom{n-j}{j}=c_{n-1,j}+\frac1{2^{n-j}}\binom{n-j}{j},
\]
for all $n\ge 2$ and $2j\le n$. Formulas \er{104} and \er{105} are used for the numerical implementations presented below. They show the dominant power-law behavior of $\vp_n$. It is possible to obtain a good approximation using the following non-rigorous reasoning. The domain of definition for $F(z)$ contains the intersection of the (filled) Julia sets related to the quadratic polynomials $rz+(1-r)z^2$ for $r\in(1/2,1)$. Thus, there is only one critical point for $F(z)$ in $\{|z|<1\}$. It is $z=1$. Substituting an assumption $F(z)\approx A(1-z)^{\a}$ into \er{101}, we obtain equation for $\a$, the analog of \er{017}:
\[\lb{106}
 1-\lt(\frac32\rt)^{\a}=-\frac38\a.
\]
The root of \er{106} with the minimal real part is $\a=-0.3904295156631794...$. This is a unique root with the negative real part. The corresponding term $A(1-z)^{\a}$ makes the greatest contribution to the growth of the coefficients $\vp_n$.  We have
\[\lb{107}
 \Phi(z)\approx-A\a(1-z)^{\a-1}=-A\a\sum_{n=0}^{+\iy}(-1)^n\binom{\a-1}{n}z^n.
\]
Using asymptotic expansion for the extended binomial coefficients \er{017}
in \er{107}, we obtain a good approximation
\[\lb{109}
 \vp_n\approx C n^{-\a},
\]
with some constant $C$, numerical approximation for which is
\[\lb{110}
 C=1.2232199....
\]
Taking more terms $F(z)\approx A(1-z)^{\a}+B(1-z)^{1+\a}$ and using the same argument as above we can improve a little bit \er{109}, namely
\[\lb{111}
 \vp_n\approx C \lt(n^{-\a}-\frac{(3\a^2+11\a+2)\a}{2(6+9\a)}n^{-\a-1}\rt).
\]
The second term $n^{-\a-1}$ is dominant in comparison with other terms related to the roots of \er{106} with positive real parts since these real parts are larger than $1+\a$. A comparison between exact values $\vp_n$, computed by \er{104} and \er{105}, and approximation terms, see \er{111}, is presented in Fig. \ref{fig1}. Two terms in \er{111} already give a good approximation of $\vp_n$ as it is seen in Figs. \ref{fig1a} and \ref{fig1b}. Double precision computations are not enough for good accuracy and lead to numerical errors appearing in Fig. \ref{fig1c}. The root $\a=0$ of \er{106} should not contribute to the asymptotics, since the corresponding $\a=-1$ is not a root of the original equation \er{017} for this example - writing \er{106} instead of \er{017}, we got a few obvious spurious roots. Thus, Fig. \ref{fig1c} should contain something attenuated instead of what we see. I recalculated Fig. \ref{fig1c} with double-double precision and with a more precise constant 
\[\lb{110dd}
 C=1.223219951386792...,
\]
see \er{110}. The result is seen in Fig. \ref{fig1cdd}. Now, it is correct. Let us only note that obtaining the most accurate values of the constants is necessary for asymptotics of this type. In the next Example 2, special attention is paid to this aspect. 

\begin{figure}
	\centering
	\begin{subfigure}{0.99\linewidth}
		\includegraphics[width=\linewidth]{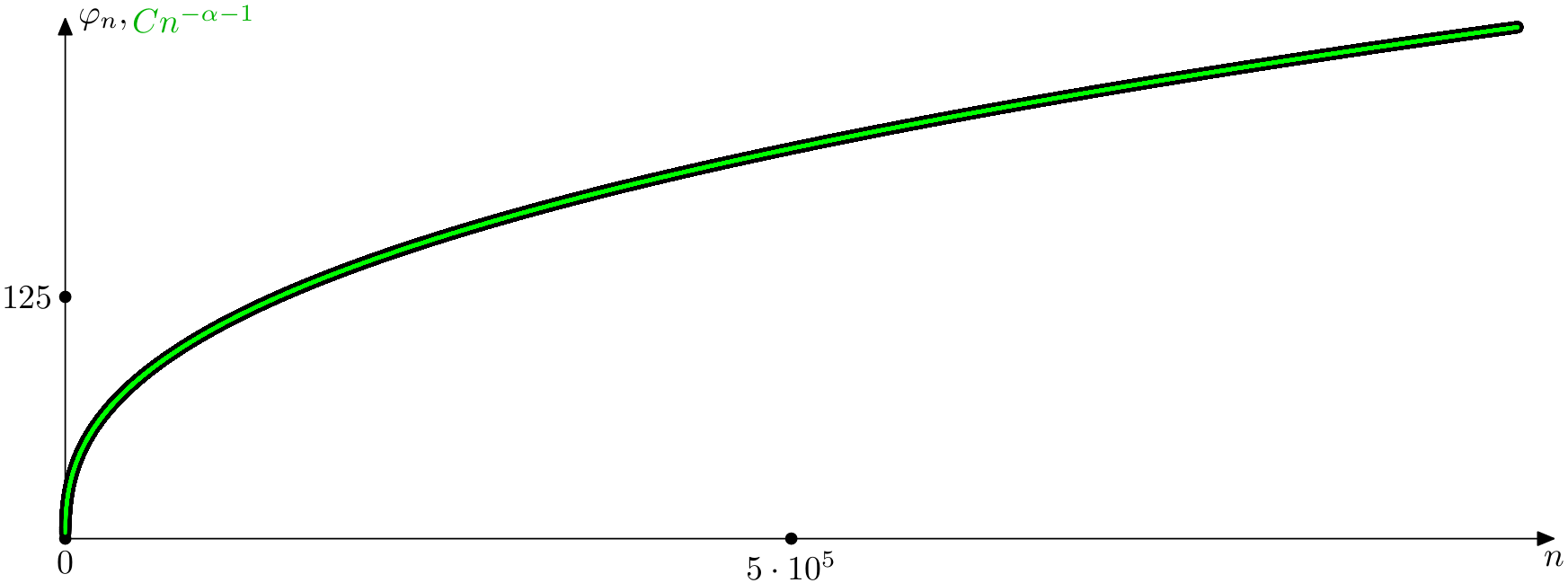}
		\caption{Exact value and first approximation term}
		\label{fig11a}
	\end{subfigure}\vfill
	\begin{subfigure}{0.99\linewidth}
		\includegraphics[width=\linewidth]{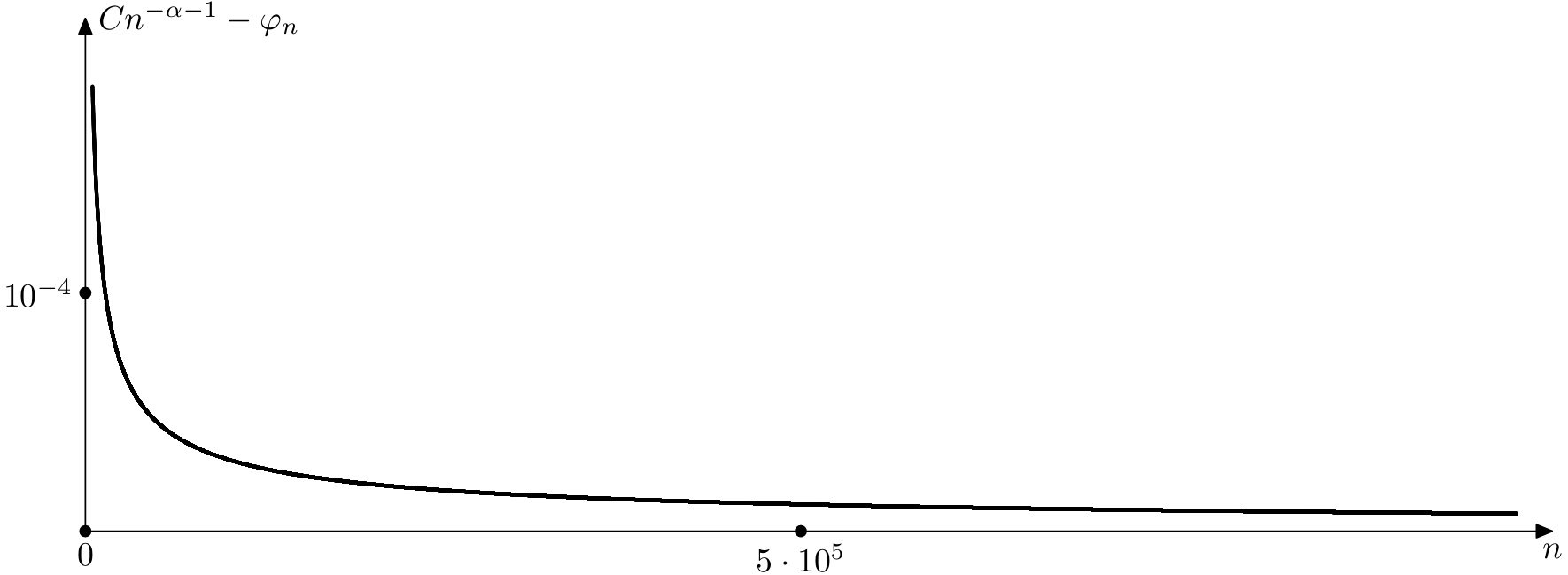}
		\caption{Difference between first approximation term and exact value}
		\label{fig11b}
	\end{subfigure}\vfill
	\begin{subfigure}{0.99\linewidth}
		\includegraphics[width=\linewidth]{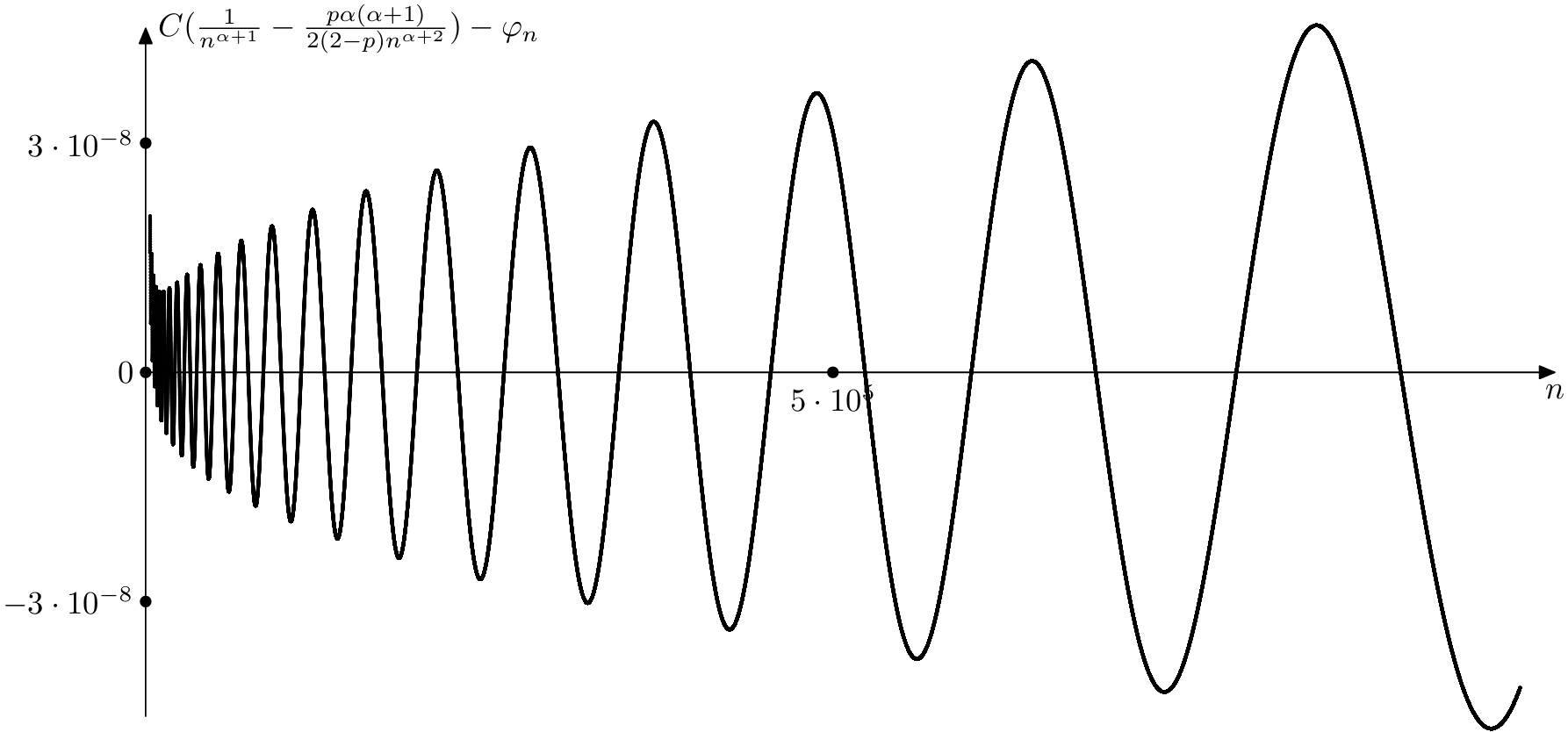}
		\caption{Difference between exact value and two first approximation terms}
		\label{fig11c}
	\end{subfigure}
	\caption{Comparison between exact values and their approximations for the fixed offspring distribution, see \er{1110} and \er{1111}.}
	\label{fig11}
\end{figure}

Before we continue, let us emulate the power law \er{109} by using one polynomial only $P(z)=pz+(1-p)z^2$ with fixed $p$. We can use the results of the last Example  3 with $p=a=b$, see \er{126}. In our case,
\[\lb{1110}
 p=0.6791281732038788538781...
\]
leads to
\[\lb{1111}
 \a=-1.3904295156631794...,\ \ \ C=1.20166998031....
\]
The power $\a$ is the same as in \er{109}, but $C$ is a little bit smaller, see \er{110dd}. However, the other asymptotic terms differ significantly in these two examples. In the example with the fixed offspring distribution, the oscillations appear even in the leading term, see Fig. \ref{fig11} and compare Fig. \ref{fig1cdd} with Fig. \ref{fig11c}, as discussed at the end of the Introduction section. Just like in the last example, we broke tradition and made all the calculations with double-double precision, bypassing double precision.

\subsection{Example 2.} The case $P_r(z)=rz+(1-r)z^2$, with uniform $r\in(0,1)$. Applying the same arguments as in the previous example, we obtain the analog of \er{101} for the integral $F(z)=\int_0^z\Phi(\z)d\z$, namely 
\[\lb{112}
\frac{F(z)-F(z^2)}{z-z^2}=\frac12F'(z),\ \ F(0)=F'(0)=0,\ \ F''(0)=1,
\]
which gives the analog of \er{103}:
\[\lb{113}
   \frac{\vp_{2n}}{2}=\frac{\vp_{2n}}{2n+1}+\frac{\vp_{2n-1}}{2n}+...+\frac{\vp_{n}}{n+1},\ \ \
   \frac{\vp_{2n+1}}{2}=\frac{\vp_{2n+1}}{2n+2}+\frac{\vp_{2n}}{2n+1}+...+\frac{\vp_{n+1}}{n+2}.
\]
RHSs in \er{113} are almost identical to each other. This fact allows us to simplify \er{113} easily
\[\lb{114}
 \vp_n=\frac{n+1}{n-1}\vp_{n-1},\ \ n\ {\rm mod}\ 2=0;\ \ \ \vp_n=\frac{n+1}{n-1}\vp_{n-1}-\frac{4}{n-1}\vp_{\frac{n-1}2},\ \ n\ {\rm mod}\ 2=1,
\]
with the initial state $\vp_1=1$. Introducing $\p_n=\vp_n/n$, we can rewrite \er{114} as
\[\lb{115}
\p_n=\lt(1+\frac1n\rt)\p_{n-1},\ \ n\ {\rm mod}\ 2=0;\ \ \ \p_n=\lt(1+\frac1n\rt)\p_{n-1}-\frac{2}{n}\p_{\frac{n-1}2},\ \ n\ {\rm mod}\ 2=1,
\]
with the initial state $\p_1=1$. By induction, one can check 
\[\lb{115a}
 \p_n=\frac{n+1}2\lt(1+\sum_k\prod_{1\le r\le k,\ 1<a_{j_r}\le n,\ 2<\frac{a_{j_r}}{a_{j_{r+1}}},\ a_{j_r}\ {\rm mod}\ 2=1}\frac{-1}{a_{j_r}}\rt).
\]
There are a few elements of this sequence
$$
 \p_2=\frac32,\ \ \p_3=\frac42\lt(1-\frac13\rt),\ \ \p_4=\frac52\lt(1-\frac13\rt),\ \ \p_5=\frac62\lt(1-\frac15-\frac13\rt),\ \ \p_6=\frac72\lt(1-\frac15-\frac13\rt),
$$
$$
 \p_7=\frac82\lt(1-\frac17-\frac15-\frac13+\frac1{7\cdot3}\rt),\ \ \p_8=\frac92\lt(1-\frac17-\frac15-\frac13+\frac1{7\cdot3}\rt),\ \ \p_{15}=\frac{16}2\lt(1-\frac1{15}-\frac1{13}-\frac1{11}-\frac19-
$$ 
$$ 
 \frac17-\frac15-\frac13+\frac1{15\cdot7}+\frac1{15\cdot5}+\frac1{13\cdot5}+\frac1{11\cdot5}+\frac1{15\cdot3}+\frac1{13\cdot3}+\frac1{11\cdot3}+\frac1{9\cdot3}+\frac1{7\cdot3}-\frac1{15\cdot7\cdot3}\rt).
$$
The analog of \er{017} and, respectively \er{106}, for \er{112} is
\[\lb{116}
 1-2^{\a}=-\frac12\a.
\]
The root of \er{116} with the minimal real part is $\a=-1$. Thus, the corresponding approximation of $F(z)\approx A(1-z)^{-1}$ leads to the linear growth of $\vp_n\approx An$, which is in good agreement with numerical computations. However, in this case $z=1$ is not a unique singularity for $F(z)$ in $\{z:\ |z|\le1\}$. One of the reasons is that the filled Julia set for $P_0(z)=z^2$ is the unit disk. Analyzing \er{112}, it is seen that $e^{\frac{2\pi\mathbf{i}k}{2^n}}$, $k,n\ge0$ are singularities for $F(z)$ giving different impacts to the power series coefficients. Let us extract the impact of two major singularities $z=1$ and $z=-1$. Substituting
\[\lb{117}
 F(z)\approx A(1-z)^{-1}+B\ln(1-z)+C\ln(1+z)
\]
into \er{112} and eliminating two main terms at $z=1$ and one main term at $z=-1$, we obtain
\[\lb{118}
 B=\frac{A}{4\ln 2-2},\ \ \ C=\frac A2.
\]
For the constant $A$, we have the limit $A=\lim_{n\to\iy}\p_n$. Researcher Tian Vla\v{s}i\'c from {\it math.stackexchange.com} notes that the knowledge system WolframAlpha shows 
\[\lb{119}
 A=\frac1{2-2\ln2},
\]
see the discussion at \href{https://math.stackexchange.com/questions/4748129/asymptotics-of-sequence-of-rational-numbers}{this site}\endnotemark[\getrefnumber{ref1}]. Identity \er{119} agrees very well with numerical data. It is useful to note that the constant $1/(1-\ln2)$ appears at least in one somewhat similar but inhomogeneous equation, see \cite{BS}. Also, equations similar to \er{114} but with constant coefficients were considered in, e.g., \cite{EHOP}. The corresponding asymptotic includes factors with the logarithm of natural numbers in the denominator.

\begin{figure}
	\centering
	\begin{subfigure}{0.99\linewidth}
		\includegraphics[width=\linewidth]{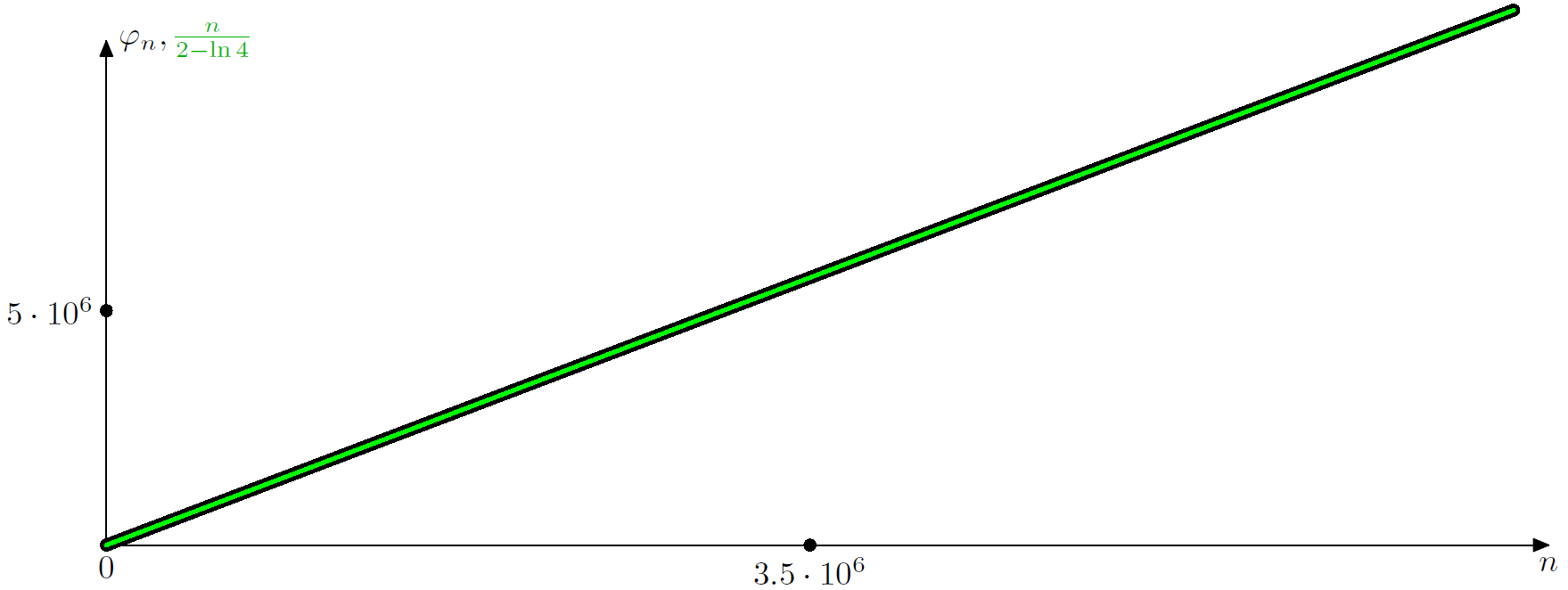}
		\caption{Exact value and first approximation term}
		\label{fig2a}
	\end{subfigure}\vfill
	\begin{subfigure}{0.99\linewidth}
		\includegraphics[width=\linewidth]{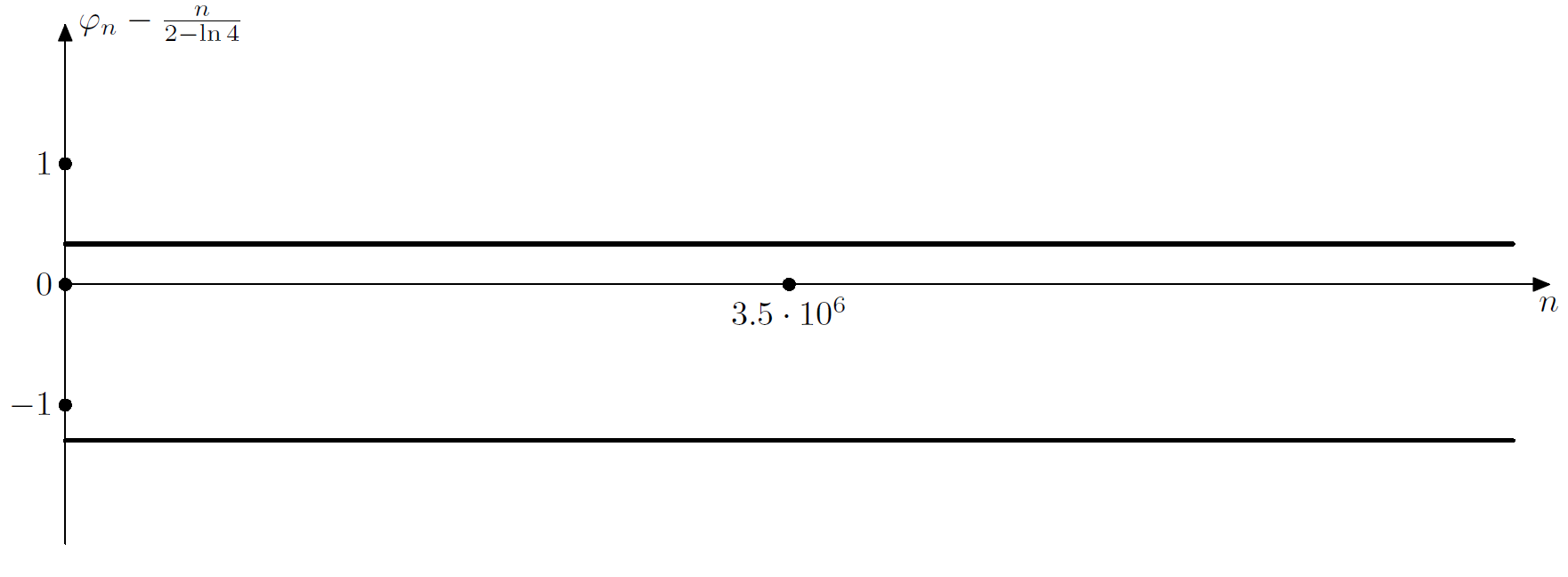}
		\caption{Difference between exact value and first approximation term}
		\label{fig2b}
	\end{subfigure}\vfill
	\begin{subfigure}{0.99\linewidth}
		\includegraphics[width=\linewidth]{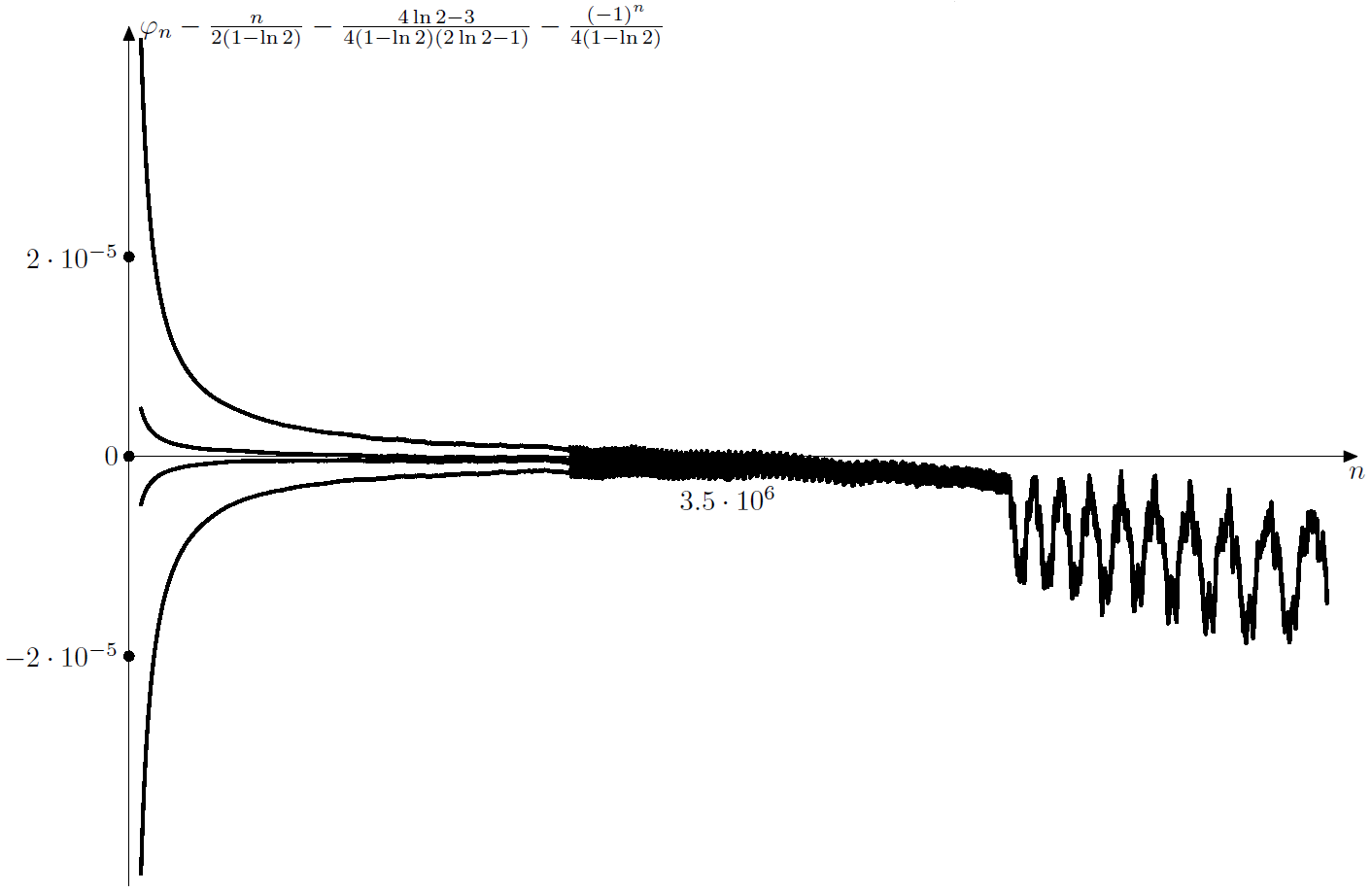}
		\caption{Difference between exact value and three approximation terms}
		\label{fig2c}
	\end{subfigure}
	\caption{Comparison between exact values and their approximations, see \er{120}.}
	\label{fig2}
\end{figure}

Finally, substituting \er{119} and \er{118} into \er{117} and using $\Phi(z)=F'(z)$, we obtain very good approximations
\[\lb{120}
 \vp_n\approx A(n+1)-B+C(-1)^n=\frac{n}{2(1-\ln2)}+\frac{4\ln2-3}{4(1-\ln2)(2\ln2-1)}+\frac{(-1)^n}{4(1-\ln2)}.
\]
A comparison between exact values $\vp_n$, computed by \er{114}, and approximation terms, see \er{120}, is presented in Fig. \ref{fig2}. Three terms in \er{120} already give a good approximation of $\vp_n$ as it is seen in Figs. \ref{fig2a} and \ref{fig2b}. Note that, the values at even and odd $n$ are different in Fig. \ref{fig2b}. Again, double precision computations are not enough for good accuracy and lead to some small numerical errors appearing in Fig. \ref{fig2c}. Nevertheless, some attenuated factors from complex singularities of $F(z)$ are seen in Fig. \ref{fig2c}. The singularities $z=\pm\mathbf{i}$ give the largest impact to the remainder. Along with the next asymptotic terms related to $z=\pm1$, they lead to oscillations of period $4$ - values of the remainder at $n$, $n+1$, $n+2$, and $n+3$ differ significantly. The corresponding exact terms can be obtained by the analogy with the use of \er{117}. On the other hand, analyzing numerical results, one can directly substitute four next asymptotic terms of the order $n^{-1}$ into \er{114}, and use the symmetry among them - two unique values and two others with the different sign - to obtain the next more accurate approximation
\[\lb{120a}
 \vp_n=\frac{1}{2(1-\ln2)}\lt(n+\frac{4\ln2-3}{4\ln2-2}+\frac{(-1)^n}{2}-\frac{(-1)^{n}(1-\ln2)}{(2\ln2-1)n}+\frac{\cos\frac{\pi n}2-\sin\frac{\pi n}2}{n}+...\rt).
\] 
It seems that the remainder in \er{120a} is of the order $n^{-2}$. Continuing the process with the next asymptotic terms, we hopefully can reach the long-phase oscillations related to the complex roots $\a$ of \er{116}. In this case, the asymptotic for $\vp_n$ will contain both: short-phase and long-phase oscillations.

\begin{figure}
	\centering
		\includegraphics[width=0.9\linewidth]{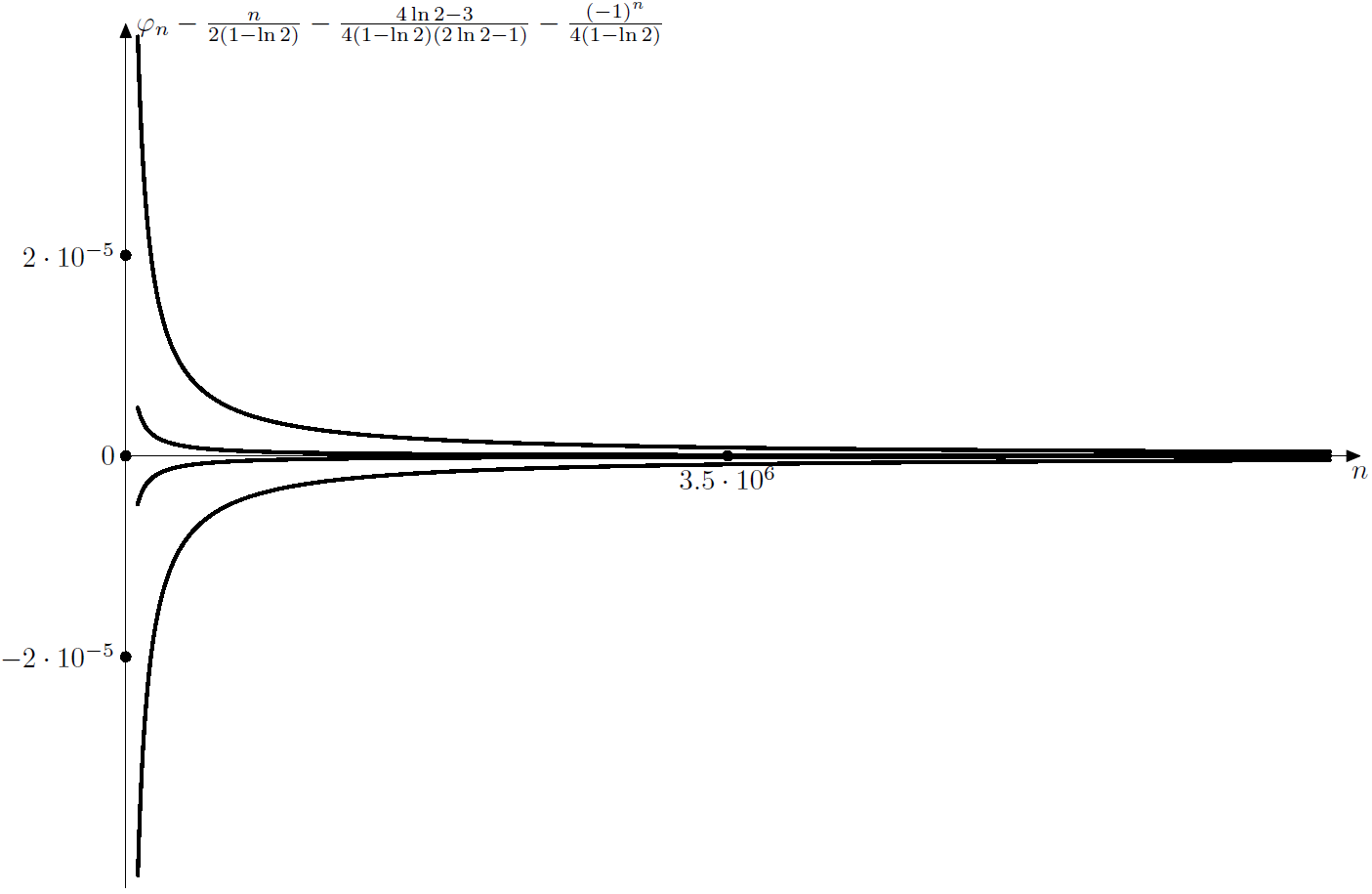}
		\caption{The same computations as in Fig. \ref{fig2c} with using double-double precision instead of double precision.}
		\label{fig2cdd}
\end{figure}

\begin{figure}
	\centering
		\includegraphics[width=0.9\linewidth]{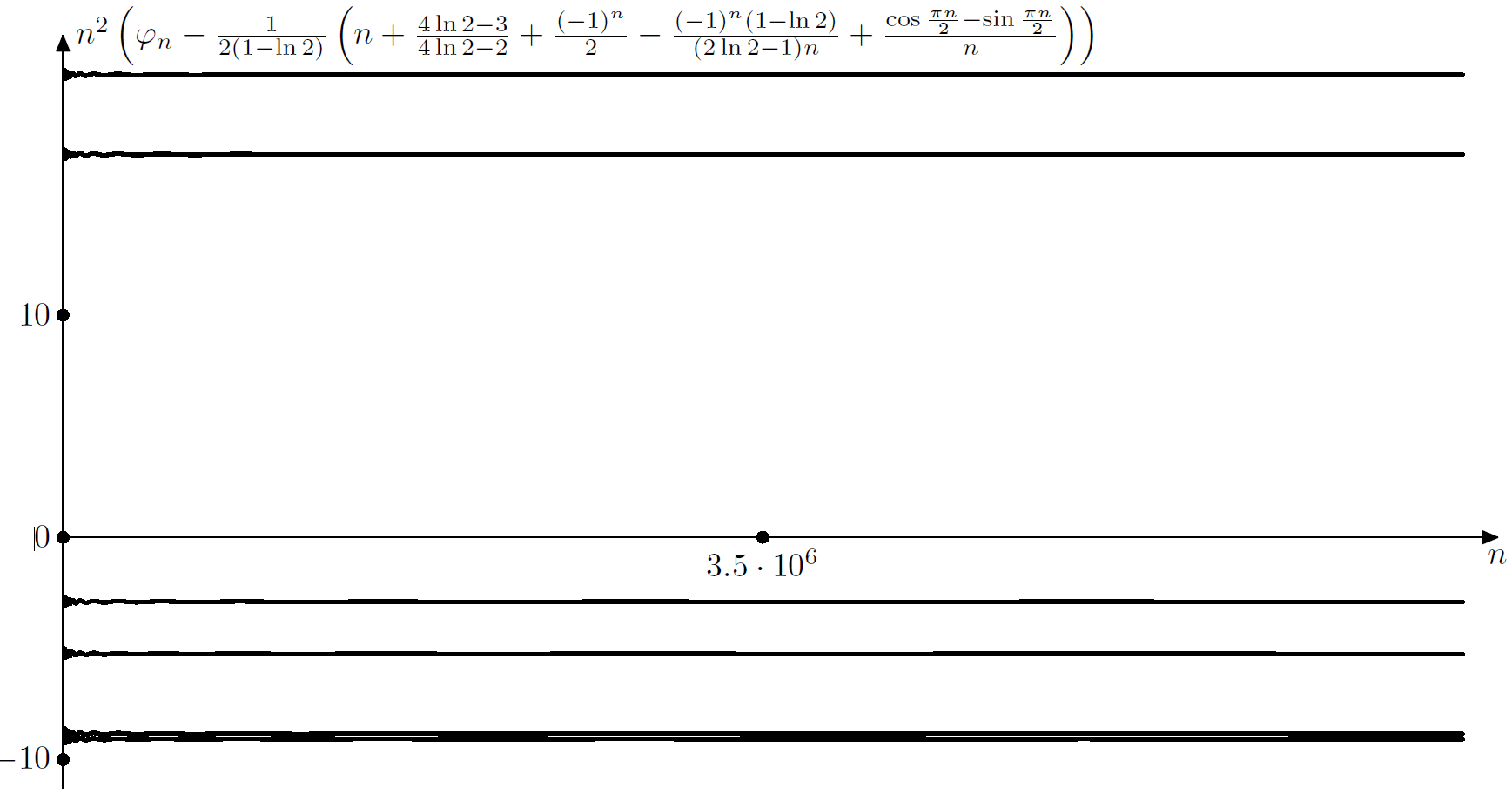}
		\caption{A few next asymptotic terms to \er{120} are added, see \er{120a}.}
		\label{fig3cdd}
\end{figure}

\begin{figure}
	\centering
	\includegraphics[width=0.85\linewidth]{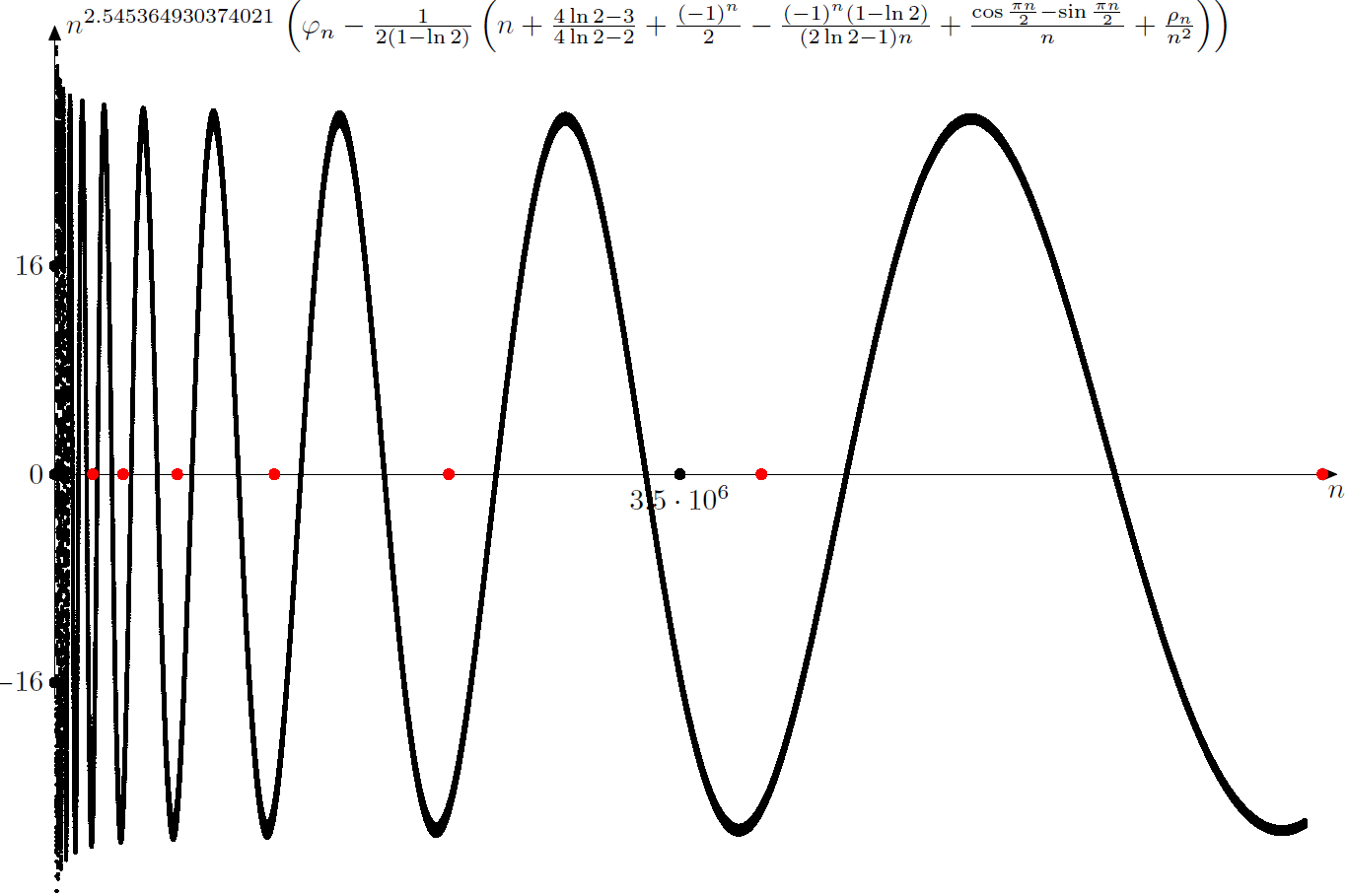}
	\caption{A few more asymptotic terms to Fig. \ref{fig3cdd} are added. Red points are $\exp(2\pi j/10.75397517526887)$ for $j=21,...,27$.}
	\label{fig4cdd}
\end{figure}

Since \er{114} is much simpler than \er{104}, we can perform computations with double-double precision in a moderate time. For the computations with a standard double precision we use \href{https://www.embarcadero.com/ru/products/delphi/starter}{Embarcadero Delphi CE} with a highly optimized parallel library \href{https://www.dewresearch.com/products/mtxvec/mtxvec-for-delphi-c-builder}{MtxVec}. For very precise computations with double-double precision, we use an additional library \href{https://github.com/neslib/Neslib.MultiPrecision}{NesLib}. While this precision is enough for the accurate numerical implementation of \er{114}, this computation is not vectorized. Anyway, Fig. \ref{fig2cdd} is an improved version of Fig. \ref{fig2c}. The double-double precision is also enough to estimate how good the next approximation \er{120a} is, see Fig. \ref{fig3cdd}. It seems that two limit values among eight in Fig. \ref{fig3cdd} coincide with some others - that is why there are only six distinct values, while the short-phase period is $8$. Moreover, the discussed long-phase oscillations are more or less visible now. Let us extract them explicitly. The values of $8$-periodic asymptotic term $\rho_n/n^2(2-2\ln2)$ multiplied by the second power of $n$ can be computed in the same way as the previous terms - by using \er{114} and the empirical fact that $\r$ has a zero average. Skipping these straightforward computations, we write the final values
$$
 \r_{8n}=\frac{11\ln2-9}{4\ln2-2},\ \ \ \r_{8n+1}=\frac{19-31\ln2}{4\ln2-2},\ \ \ \r_{8n+2}=\frac{27-45\ln2}{4\ln2-2},\ \ \ \r_{8n+3}=\frac{\ln2-5}{4\ln2-2},
$$ 
$$
\r_{8n+4}=\frac{11\ln2-9}{4\ln2-2},\ \ \ \r_{8n+5}=\frac{33\ln2-13}{4\ln2-2},\ \ \ \r_{8n+6}=\frac{19\ln2-5}{4\ln2-2},\ \ \ \r_{8n+7}=\frac{\ln2-5}{4\ln2-2}.
$$
In principle, we can continue to compute the next asymptotic terms related to the attenuated factors $n^{-j}$, where $j>2$ is an integer. However, there are complex zeros of \er{116}
$$
 \a=2.545364930374021...\pm10.75397517526887...\mathbf{i}
$$ 
with a non-integer real part less than $3$. This is the smallest real part of zeros lying in the right half plane $\Re\a>0$. According to the discussion around \er{017} and \er{018}, these zeros contribute to a long-phase oscillation. A double-double precision still allows us to compute this contribution explicitly, see Fig. \ref{fig4cdd} 

\begin{figure}
    \centering
    \begin{subfigure}[b]{0.45\textwidth}
        \includegraphics[width=\textwidth]{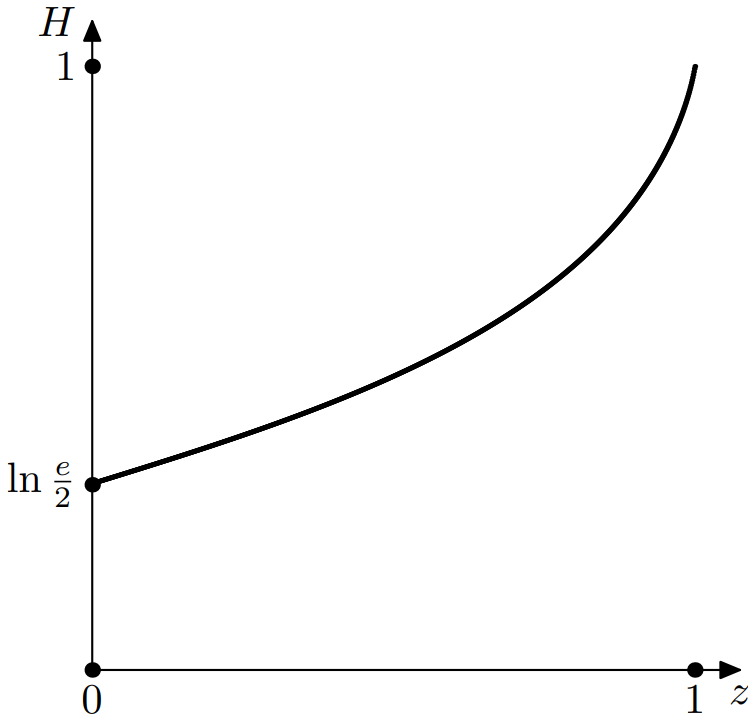}
        \caption{Initial function $H_0\ev1$ is constant in the whole interval $[0,1]$.}
        \label{fig5iter1}
    \end{subfigure}
    ~ 
    \begin{subfigure}[b]{0.45\textwidth}
        \includegraphics[width=\textwidth]{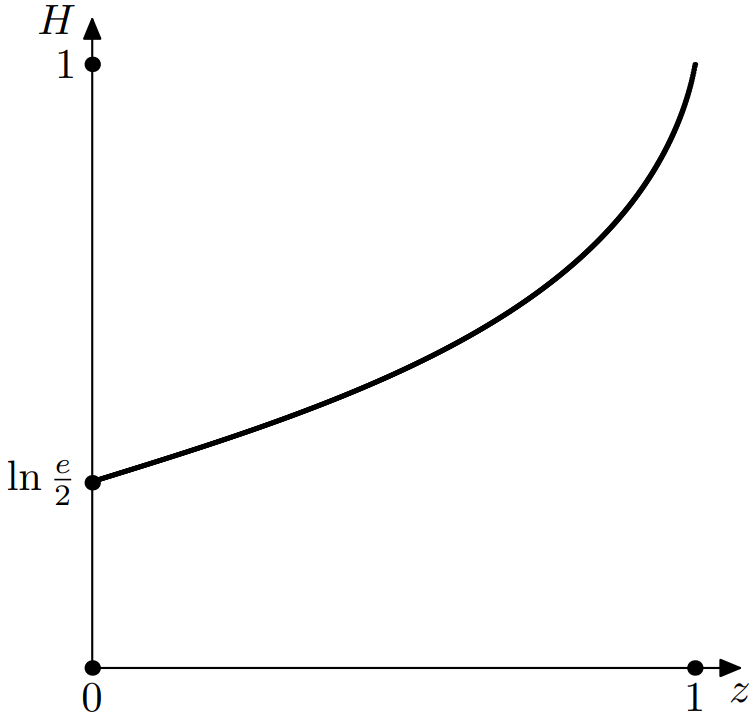}
        \caption{Initial function $H_0(z)=0$ for $z\in[0,0.5)$, and $H_0(z)=2-z$ for $z\in(0.5,1]$.}
        \label{fig5iter2}
    \end{subfigure}
    \caption{Result of $100$ Picard-Lindel\"of iterations of \er{n004} with the integration step $\D z=10^{-4}$ (rectangle method) and different starting functions $H_0$.}\label{fig5iter}
\end{figure}

Let us make a few remarks regarding the computation of $F(z)$. For the new function
\[\lb{n000}
 H(z)=\frac{(1-z)F(z)}{z^2},
\]
equation \er{112} can be rewritten in the integral form
\[\lb{n001}
 H(z)=\frac1{2(1-z)}-\frac2{1-z}\int_0^z\frac{\zeta}{1+\zeta} H(\zeta^2)d\zeta,\ \ \ H(0)=\frac12.
\]
For $|z|<1$, since $|z^2|<|z|$ (when $z\ne0$), one can use the Picard-Lindel\"of iterations
\[\lb{n002}
 H_0\equiv\frac12,\ \ \ H_{n+1}(z)=\frac1{2(1-z)}-\frac2{1-z}\int_0^z\frac{\zeta }{1+\zeta}H_n(\zeta^2)d\zeta
\]
to obtain the converging approximations $H_n(z)\to H(z)$. It also gives some new sequences converging to $H(1)$, for which we believe $H(1)=1/2(1-\ln2)$. On the other hand, if $H(1)$ is bounded as it follows from \er{117} and \er{n000}, then \er{n001} leads to
\[\lb{n003}
 \int_0^1\frac{\zeta }{1+\zeta}H(\zeta^2)d\zeta=\frac14,
\]
and, hence,
\[\lb{n004}
 H(z)=\frac1{1-z}\int_z^1\frac{2\zeta }{1+\zeta}H(\zeta^2)d\zeta.
\]
Again, starting with $H_0\equiv H(1)$ in \er{n004}, the Picard-Lindel\"of iterations give converging approximations for $H(z)$, and some new sequence for $H(1)$ if we take into account the fact that $H(0)=\frac12$. In other words, is it true that
\[\lb{n005}
 \lim_{n\to+\iy}\int_0^1\int_{z_n^2}^1...\int_{z_3^2}^1\int_{z_2^2}^1\frac{2z_1}{1+z_1}\prod_{k=2}^n\frac{2z_k}{(1+z_k)(1-z_k^2)}dz_1dz_2...dz_n=1-\ln2?
\]
This is an equivalent problem to \er{119}. In fact, for the iterations, instead of $H_0\ev1$, one can use any function such that $H_0(1)=1$ and $H_0(z)$ is continuous at $z=1$. As an example, the result of iterations for two different $H_0$ is given in Fig. \ref{fig5iter}. There are no visible differences between Figs. \ref{fig5iter1} and \ref{fig5iter2}. Hence, there is an idea that reformulates the problem about the constant \er{119}: try to find a function $H_0(z)$ continuous at $z=1$ such that $H_0(1)\ne0$ and integrals
\[\lb{n006}
 ...\int_{z_4^2}^1\frac{2z_3 }{1+z_3}\cdot\frac1{1-z_3^2}\int_{z_3^2}^1\frac{2z_2 }{1+z_2}\cdot\frac1{1-z_2^2}\int_{z_2^2}^1\frac{2z_1 }{1+z_1}H_0(z_1^2)dz_1dz_2dz_3...
\]
can be computed explicitly. 

In the end, the famous researcher Fedja found a solution of \er{119} by using a very spectacular method, see \href{https://mathoverflow.net/questions/458885/simple-integral-equation}{this site}\endnotemark[\getrefnumber{ref2}]. After reading his Proof, I find my own very different and simple Proof based on the integration by parts of \er{n004}
\begin{multline}\lb{alt1}
 \int_0^{1-\ve}\frac{H(z)}{1-z}dz=\int_0^{1-\ve}\frac1{(1-z)^2}\int_z^1\frac{2\zeta }{1+\zeta}H(\zeta^2)d\zeta dz=
\\
 \frac1{1-z}\int_z^1\frac{2\zeta }{1+\zeta}H(\zeta^2)d\zeta\bigg\rvert_{z=0}^{z=1-\ve}+\int_0^{1-\ve}\frac{2\zeta}{1-\zeta^2}H(\zeta^2)d\zeta,
\end{multline}
which leads to
\[\lb{alt2}
 \int_{(1-\ve)^2}^{1-\ve}\frac{H(z)}{1-z}dz= H(1-\ve)-H(0).
\]
Using continuity of $H(z)$ at $z=1$ and $\int_{(1-\ve)^2}^{1-\ve}(1-z)^{-1}dz\to\ln2$ for $\ve\to0$ we obtain the result.

\subsection{Example 3.} 
\begin{figure}
	\centering
	\begin{subfigure}{0.99\linewidth}
		\includegraphics[width=\linewidth]{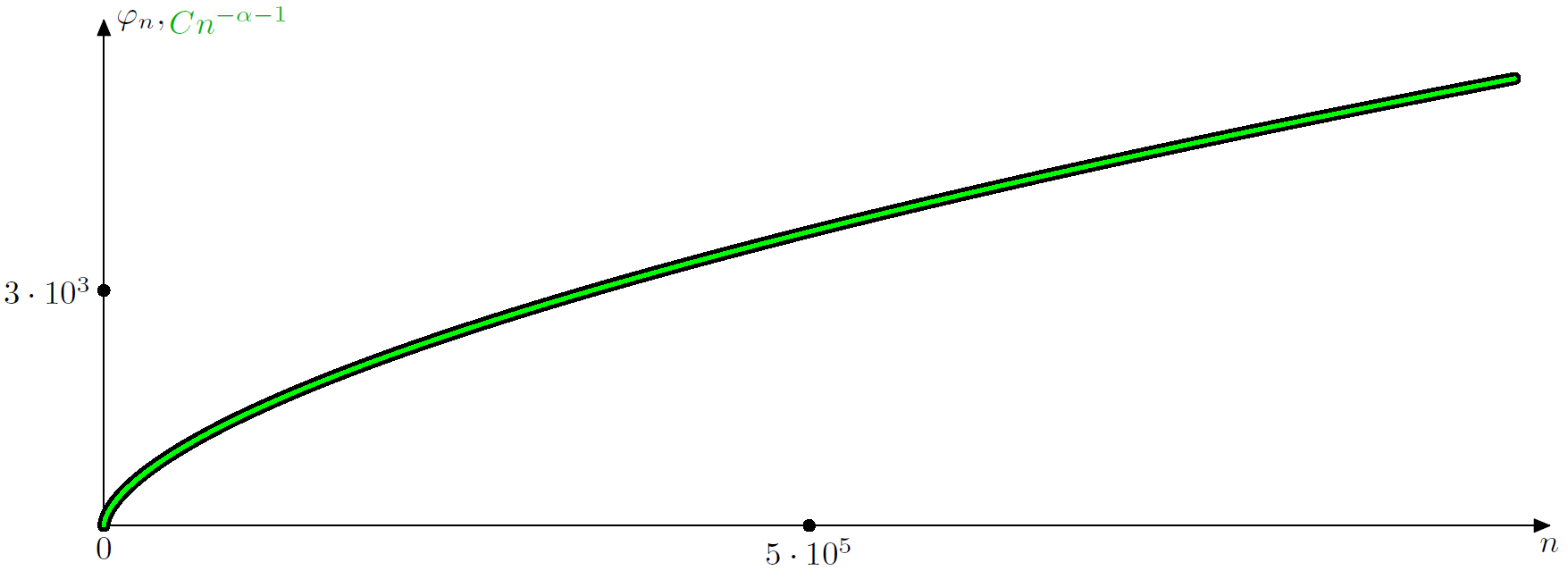}
		\caption{Exact value and first approximation term}
		\label{fig3a}
	\end{subfigure}\vfill
	\begin{subfigure}{0.99\linewidth}
		\includegraphics[width=\linewidth]{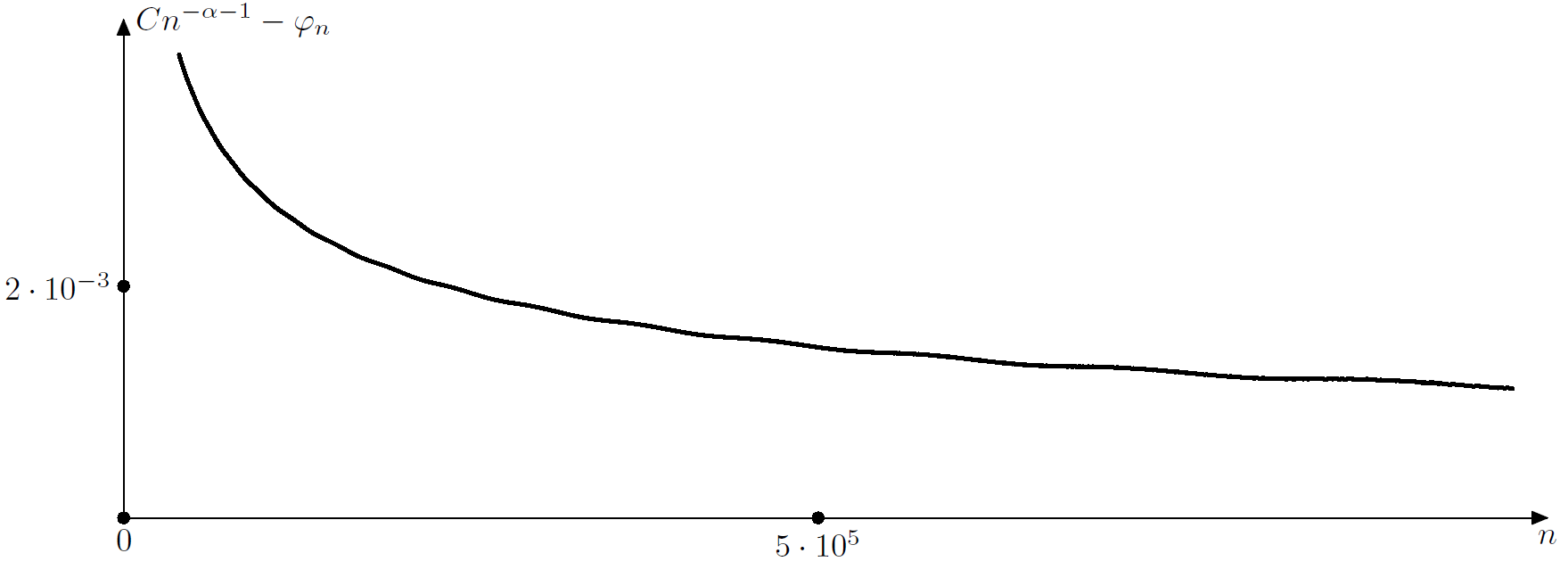}
		\caption{Difference between exact value and first approximation term}
		\label{fig3b}
	\end{subfigure}\vfill
	\begin{subfigure}{0.99\linewidth}
		\includegraphics[width=\linewidth]{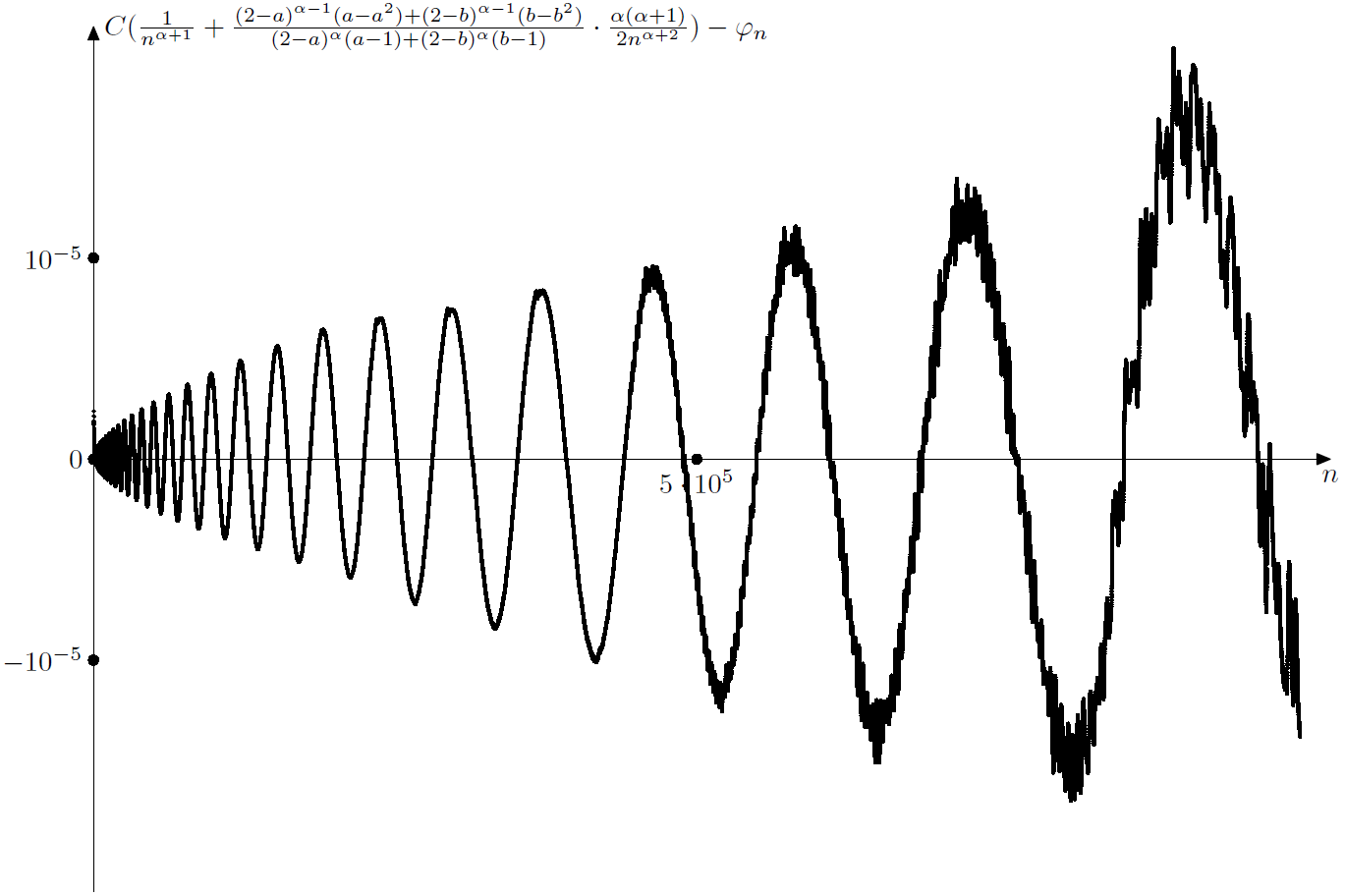}
		\caption{Difference between exact value and two approximation terms}
		\label{fig3c}
	\end{subfigure}
	\caption{Comparison between exact values and their approximations for the case \er{126}-\er{129}.}
	\label{fig3}
\end{figure}

\begin{figure}
	\centering
	\includegraphics[width=0.99\linewidth]{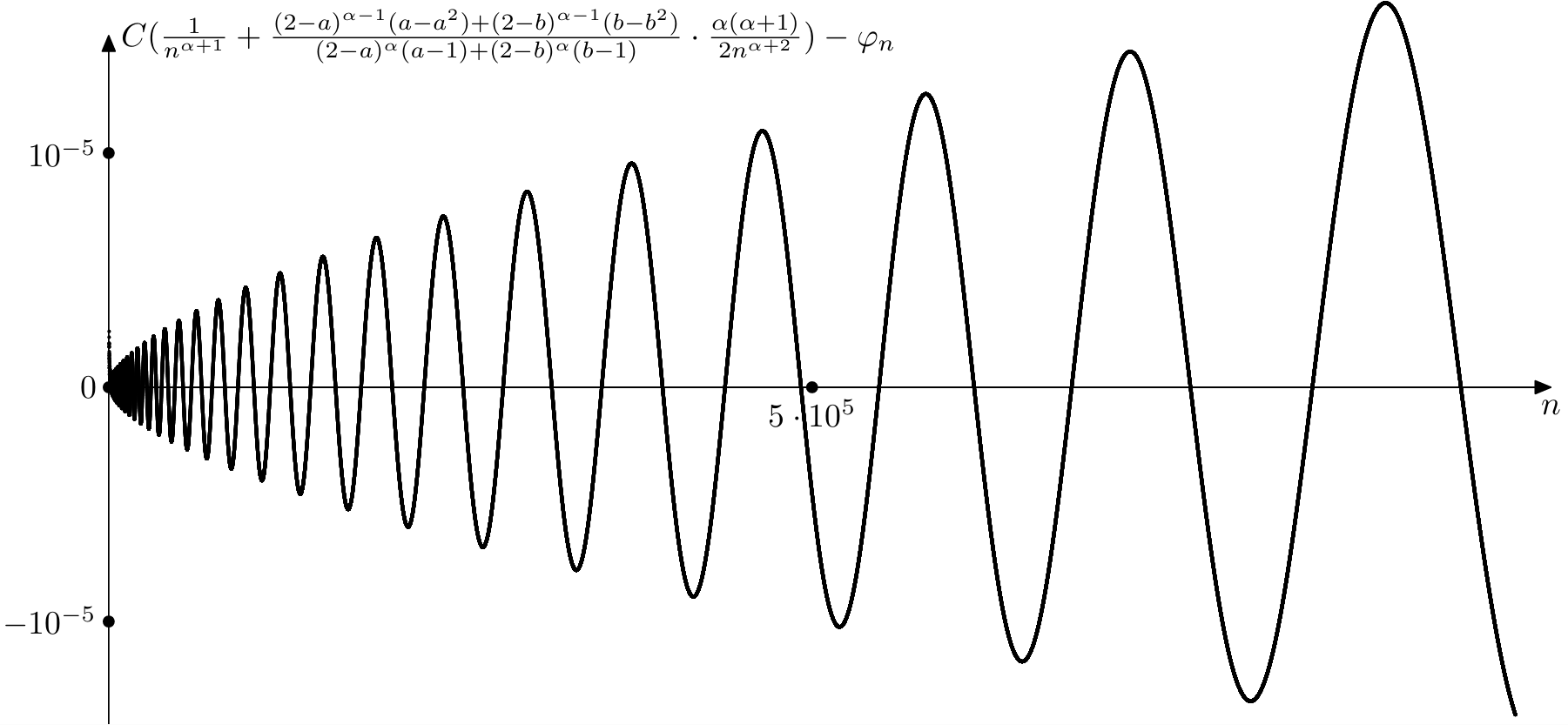}
	\caption{Correction of Fig. \ref{fig3c} with double-double precision.}
	\label{fig3cdd1}
\end{figure}

Let us consider the case of two polynomials $P_1(z)=az+(1-a)z^2$ and $P_2(z)=bz+(1-b)z^2$, $a,b\in(0,1)$, which appear with equal probability at each time step. Equation \er{007} becomes
\[\lb{121}
 \Phi(az+(1-a)z^2)+\Phi(bz+(1-b)z^2)=(a+b)\Phi(z).
\]
By analogy with the previous examples, we obtain the explicit recurrence expression for $\vp_n$:
\begin{multline}\lb{121a}
 (a+b-a^n-b^n)\vp_n=(a^{n-2}(1-a)^1+b^{n-2}(1-b)^1)\binom{n-1}{1}\vp_{n-1}+\\
 (a^{n-4}(1-a)^2+b^{n-4}(1-b)^2)\binom{n-2}{2}\vp_{n-2}+(a^{n-6}(1-a)^3+b^{n-6}(1-b)^3)\binom{n-3}{3}\vp_{n-3}+...,
\end{multline}
where the sum breaks when the lower index becomes strictly greater than the upper index in the binomial coefficient. We substitute an {\it ansatz}
\[\lb{122}
 \Phi(z)\approx A(1-z)^{\a}+B(1-z)^{\a+1}
\]
into \er{121} to obtain the analog of \er{015} for $\a$:
\[\lb{123}
 (2-a)^{\a}+(2-b)^{\a}=a+b,
\]
and another one connecting $A$ and $B$:
\[\lb{124}
 A(2-a)^{\a}(-\a)\frac{1-a}{2-a}+B(2-a)^{\a+1}+A(2-b)^{\a}(-\a)\frac{1-b}{2-b}+B(2-b)^{\a+1}=B(a+b),
\]
which leads to
\[\lb{125}
 B=A\frac{-\a(2-a)^{\a-1}(1-a)-\a(2-b)^{\a-1}(1-b)}{a+b-(2-a)^{\a+1}-(2-b)^{\a+1}}.
\]
Taking into account arguments \er{016}, \er{017}, along with \er{122}, \er{123}, and \er{125}, we obtain
\[\lb{126}
 \vp_n\approx C\lt(\frac1{n^{\a+1}}+\frac{(2-a)^{\a-1}(a-a^2)+(2-b)^{\a-1}(b-b^2)}{(2-a)^{\a}(a-1)+(2-b)^{\a}(b-1)}\cdot\frac{\a(\a+1)}{2n^{\a+2}}\rt),
\]
where $C=A/\Gamma(-\a)$. Now, the concrete values
\[\lb{127}
 a=\frac{7}{16},\ \ \ b=\frac{3}{4}
\]
are chosen such that \er{123} can be solved explicitly
\[\lb{128}
 \a=\log_{\frac54}\frac{-1\pm\sqrt{\frac{23}{4}}}{2}.
\]
For the parameter $C$ in \er{126}, numerical computations give
\[\lb{129}
 C\approx1.3355247475.
\]
There is one real $\a$ in \er{128} with the minimal real part, chosen for the approximation \er{126}, and infinitely many complex $\a$ with the same minimal real part, and other complex values. As it is seen in Fig. \ref{fig3}, two terms in \er{126} give a good approximation of $\vp_n$. Again, double precision does not provide accurate results in Fig. \ref{fig3c}, but acceptable. The corresponding double-double precision corrections are given in Fig. \ref{fig3cdd1}  The long-phase logarithmic oscillations coming from the complex $\a$, discussed above, will be dominant in comparison with the second term for very large $n$. This completely resembles the situation with the classical Galton--Watson process described in, e.g., \cite{K}. 
The difference is that in the classical case of one polynomial, we have a complete asymptotic series for $\vp_n$ and corresponding fast algorithms for the computation of the asymptotic terms. In the considered mixed case of two polynomials, there are no direct analogs of Karlin--McGregor functions and, at the moment, I do not know fast algorithms for the computation of coefficients in \er{019} - only in some cases and for some of them it is sometimes possible to give explicit formulas, as in the examples above.

\begin{figure}
	\centering
	\begin{subfigure}{0.99\linewidth}
		\includegraphics[width=\linewidth]{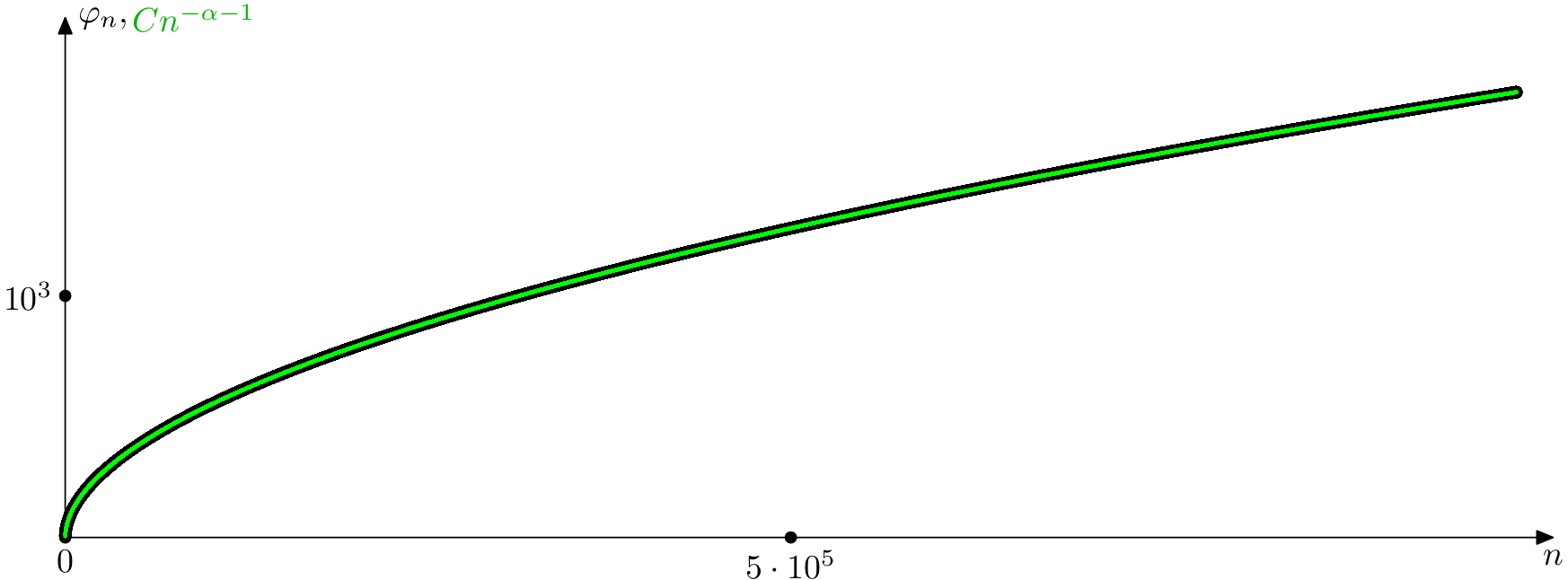}
		\caption{Exact value and first approximation term}
		\label{fig4a}
	\end{subfigure}\vfill
	\begin{subfigure}{0.99\linewidth}
		\includegraphics[width=\linewidth]{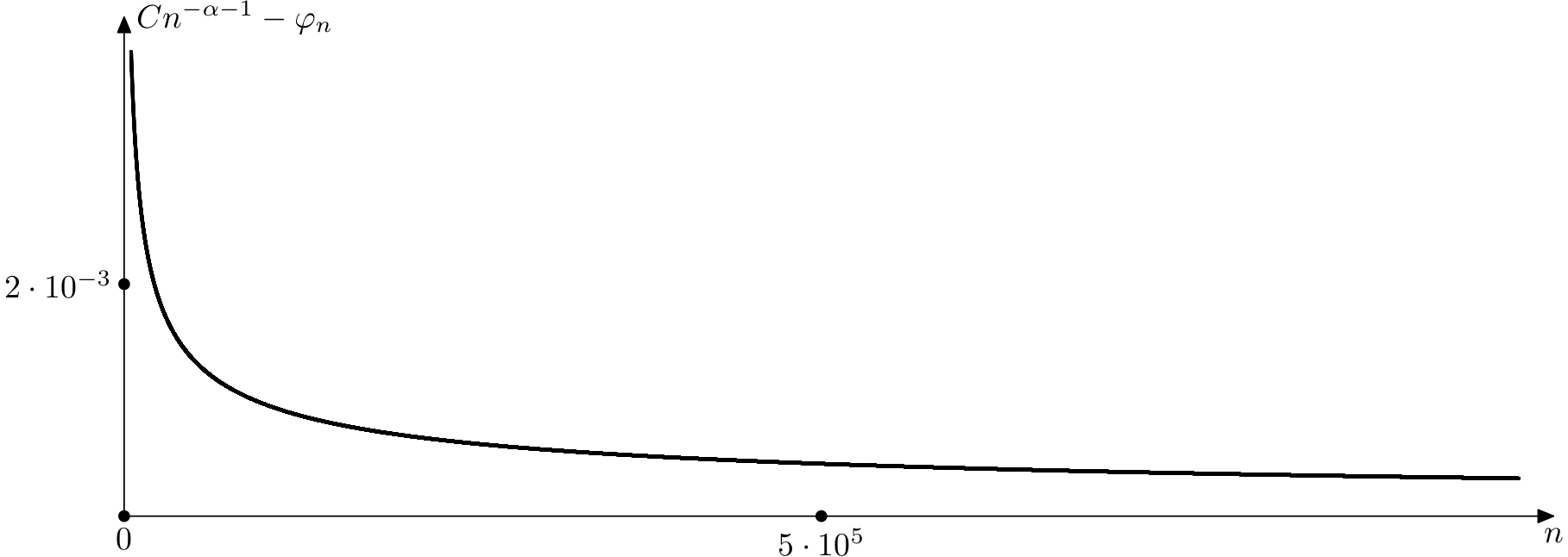}
		\caption{Difference between exact value and first approximation term}
		\label{fig4b}
	\end{subfigure}\vfill
	\begin{subfigure}{0.99\linewidth}
		\includegraphics[width=\linewidth]{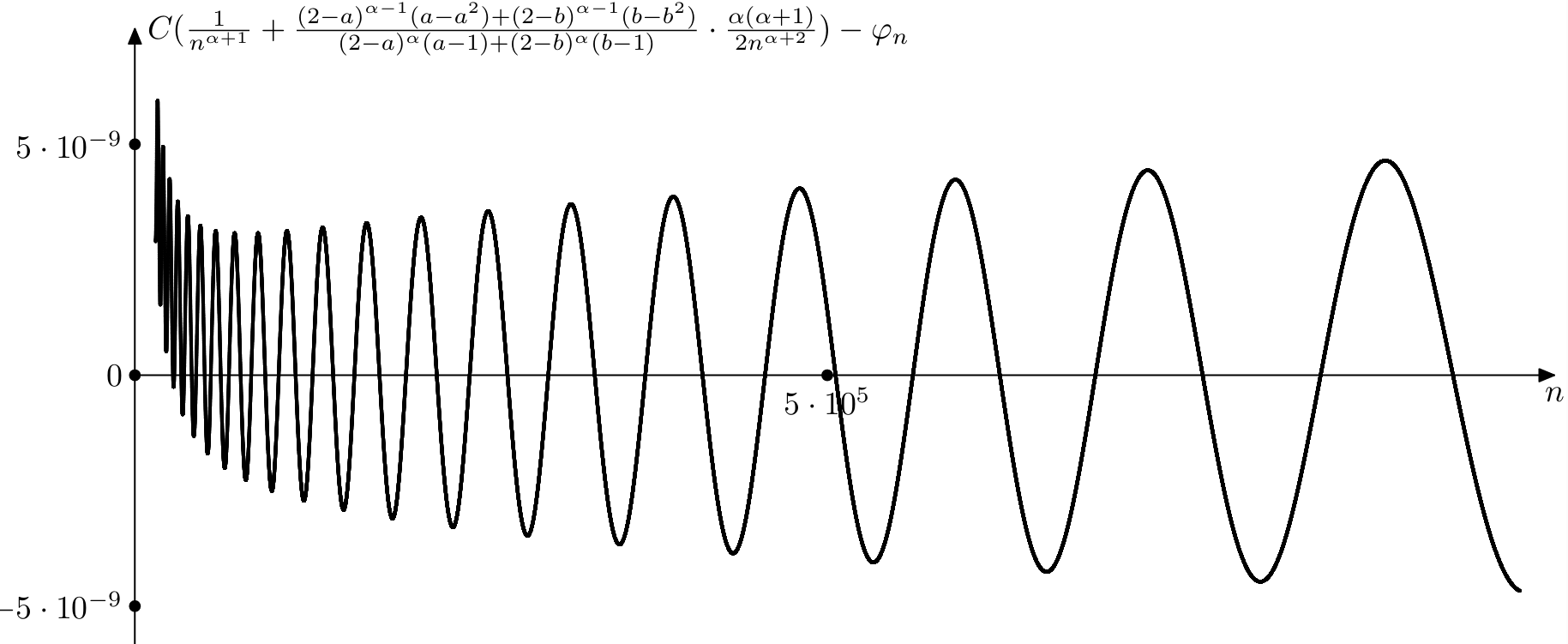}
		\caption{Difference between exact value and two approximation terms}
		\label{fig4c}
	\end{subfigure}
	\caption{Comparison between exact values and their approximations for the case \er{130}-\er{132}.}
	\label{fig4}
\end{figure}

Now, let us change a little bit the probability in the first polynomial and take 
\[\lb{130}
 a=\frac{1}{2},\ \ \ b=\frac{3}{4},
\] 
see \er{127}. In this case, \er{123} cannot be solved explicitly. It has one zero with a minimal real part
\[\lb{131}
 \a=-1.526066812384411....
\]
The corresponding factor in the main asymptotic term $\vp_n=Cn^{-\a-1}+...$ is approximately 
\[\lb{132}
 C\approx 1.28574621970439
\]
All other zeros of \er{123} are complex and lie in a strip. Real parts of some of the zeros are arbitrarily close to $\a$ mentioned in \er{131}. They make a big contribution to the asymptotic of $\vp_n$, but for very large $n$, since their factors $C_{\a}$ are tiny. For the moderate values of $n$, other $\a$ contribute more. We plot the corresponding approximations in Fig. \ref{fig4}. Breaking tradition, we performed all the calculations with double-double precision. In contrast to Fig. \ref{fig3c}, the oscillation in Fig. \ref{fig4c} has a smaller order of amplification than the main asymptotic term presented in Fig. \ref{fig4a}. In addition, if you are interested in learning about the distribution of zeros in equations of type \er{015} with a finite discrete $(R,\m)$, you can see, e.g., \cite{HW} and references therein.

\section{Conclusion}

Galton--Watson processes in random environments, when probabilities of fission at each time step are random variables themselves, are considered. The corresponding generalized Schr\"oder type functional equation describing the relative limit densities of descendants $\vp_n$ is derived. For the classical Galton--Watson process with fixed probabilities of fission, the corresponding asymptotic is $\vp_n=\sum_{\a,j}C_{\a,j}n^{-\a-j}$, where all the values $\a$ have the same real part and linearly distributed imaginary parts. Moreover, there are fast numerical procedures for the computation of $C_{\a,j}$. In contrast to the classical case, in the generalized asymptotic the real parts of $\a$ may differ, and the distribution of both imaginary and real parts can be irregular. The corresponding examples are considered. At the moment, I do not know fast algorithms for the computation of $C_{\a,j}$ in the generalized cases. These cases are the sources of many interesting analytic problems related to special functions and constants. One of them is considered in Example 2 announced also on \href{https://math.stackexchange.com/questions/4748129/asymptotics-of-sequence-of-rational-numbers}{this site}\endnotemark[\getrefnumber{ref1}].  

\section*{Acknowledgements} 
This paper is a contribution to the project M3 of the Collaborative Research Centre TRR 181 "Energy Transfer in Atmosphere and Ocean" funded by the Deutsche Forschungsgemeinschaft (DFG, German Research Foundation) - Projektnummer 274762653. 

I would like to thank users of {\it math.stackexchange.com} and {\it mathoverflow.net} for very useful discussions. 

%

\theendnotes

\end{document}